\pdfoutput=1
\documentclass[onecolumn,notitlepage,11pt]{article}

\usepackage[utf8]{inputenc}
\usepackage[T1]{fontenc}
\usepackage{ae}
\usepackage{aecompl}
\usepackage[english]{babel}
\usepackage{csquotes}
\usepackage{sectsty}
\usepackage{mathrsfs}
\usepackage{latexsym}
\usepackage{amsbsy}
\usepackage{amssymb}
\usepackage{verbatim}
\usepackage{amsthm}
\usepackage{amsfonts}
\usepackage{amsmath}
\usepackage{thmtools}
\usepackage{layout}
\usepackage{pb-diagram}
\usepackage{amscd}
\usepackage{anysize}
\usepackage{tikz}
\usepackage{graphicx}
\usepackage[all,cmtip]{xy}
\usepackage{youngtab}
\usepackage{enumerate}
\usepackage[protrusion=true,expansion=true]{microtype}
\usepackage{tocloft}

\newcommand{\dd}{\textup{d}} 
 
\newcommand{\C}{\mathbb{C}} 
\newcommand{\R}{\mathbb{R}} 
\newcommand{\Q}{\mathbb{Q}} 
\newcommand{\Z}{\mathbb{Z}}
\newcommand{\K}{\mathbb{K}}
\newcommand{\N}{\mathbb{N}}

\newcommand{\E}{\mathscr{E}}
\newcommand{\h}{\mathscr{H}}
\newcommand{\G}{\mathscr{G}}
\newcommand{\ch}{\mathrm{ch}}
\newcommand{\Ch}{\mathrm{Ch}}
\newcommand{\CS}{\mathrm{CS}}
\newcommand{\cs}{\mathrm{cs}}
\newcommand{\id}{\mathrm{id}}

\newcommand{\ev}{\mathrm{ev}}

\newcommand{\glres}{\mathrm{GL}_{\mathrm{res}}}
\newcommand{\ures}{\mathrm{U}_{\mathrm{res}}}
\newcommand{\lures}{\mathfrak{u}_{\mathrm{res}}}
\newcommand{\lglres}{\mathfrak{gl}_{\mathrm{res}}}
\newcommand{\grres}{\mathrm{Gr}_{\mathrm{res}}}
\newcommand{\stres}{\mathrm{St}_{\mathrm{res}}}

\newcommand{\fred}{\mathrm{Fred}}
\newcommand{\flip}{\mathrm{flip}}

\newcommand{\tr}{\mathrm{tr}}

\newcommand\restr[2]{{%
  \left.\kern-\nulldelimiterspace %
  #1 %
  \vphantom{\big|} %
  \right|_{#2} %
  }}

\newtheorem{thm}{Theorem}[section]

\declaretheorem[style=definition,name=Definition,numberwithin=section]{dfn}

\declaretheorem[style=plain,sibling=dfn,name=Lemma]{lem}
\declaretheorem[style=plain,sibling=dfn,name=Corollary]{cor}
\declaretheorem[style=remark,sibling=dfn,name=Remark]{rem}
\declaretheorem[style=definition,sibling=dfn,name=Proposition]{prop}

\usepackage[fixamsmath,disallowspaces]{mathtools}

\usepackage{todonotes}

\usepackage{tikz-cd}
\usepackage[ 
  style=alphabetic, %
  natbib=true, %
  backend=biber,   %
  bibwarn=true,     %
  texencoding=auto, %
  bibencoding=auto, %
  giveninits=true
]{biblatex}   

\renewbibmacro*{in:}{%
  \setunit{\addcomma\space}%
  \printtext{%
    \intitlepunct}}

\DeclareFieldFormat{pages}{#1}
\title{Smooth classifying spaces for differential $K$-theory}
\author{Eric Schlarmann}
\newcommand{\Addresses}{{%
\bigskip
\footnotesize
\noindent\textsc{Eric Schlarmann}\par\nopagebreak\textsc{Institut f\"ur
Mathematik, Universit\"at Augsburg}\par\nopagebreak
\texttt{EMAIL: eric.schlarmann@math.uni-augsburg.de}
}}

\date{}

\DeclareSourcemap{
  \maps[datatype=bibtex]{
    \map[overwrite]{
      \step[fieldsource=doi, final]
      \step[fieldset=url, null]
      \step[fieldset=eprint, null]
      \step[fieldset=urldate, null]
      \step[fieldset=isbn, null]
      \step[fieldset=issn, null]
      \step[fieldset=doi, null]
    }  
  }
}
\DeclareMathOperator*{\colim}{colim}

\makeatletter
\def\blfootnote{\gdef\@thefnmark{}\@footnotetext}
\makeatother

\begin{document}
\maketitle

\blfootnote{\textup{2010} \textit{Mathematics Subject Classification}: 19L50
(Primary) 19L10, 55N15, 47B10 (Secondary)}
\abstract{
    We construct a version of differential $K$-theory based on smooth Banach
    manifold models for the homotopy types $B \mathrm U\times \Z$ and $\mathrm
    U$ that appear in the topological $K$-theory spectrum. These manifolds carry
    natural differential forms that
    refine the topological universal Chern character, together with natural addition and
    inversion operations that induce the respective structure on $\hat K$.
    Our models are norm completions of the usual stable
    Grassmannian and the stable unitary
    group. Their regularity allows us to work completely on the level of
    classifying spaces, and therefore we do not need a compactness assumption on
    our manifolds that is present in many other descriptions. The constructed groups
    $\hat K(M)$ are isomorphic to
    the unique differential extension of $K$-theory that admits an $S^1$-integration.}

\tableofcontents

\section{Introduction}

For a given cohomology theory $E$ restricted to the category of
smooth manifolds, a differential refinement $\hat E$ provides a theory which
makes use of the additional geometric information. In the case of
topological $K$-theory, if a cycle is given by a vector
bundle, then a lift to a class in $\hat K$ would be defined by the
additional data of a connection. This connection will refine the Chern
character of this bundle, normally only well-defined as a cohomology class, to a
differential form. There is a set of axioms analogous to the Eilenberg--Steenrod
axioms for cohomology that characterizes such extensions, given by Bunke and
Schick \cite[Def. 1.1]{bunke2010uniqueness}. For any smooth manifold $M$, we have a
diagram of abelian groups
\[\begin{tikzcd}
        \hat K^*(M) \arrow[d, "R"]\arrow[r, "I"] & K^*(M)\arrow[d, "\ch"] \\ 
        \Omega^*_{d=0}(M) \arrow[r, "\mathrm{Rham}"]& H^*(M),
\end{tikzcd}\]
where $I$ and $R$ are forgetful functors that must come with any definition of
$\hat K$, $\ch$ is the topological Chern
character, and $\mathrm{Rham}$ is the de Rham map. Although this is not a
cartesian diagram, the slogan still is that we combine $K$-theory and forms in a
(homotopy theoretic) fiber product
\begin{align*}
    "\mathrm{Differential}\; K\text{-theory} = K\text{-theory}
    \times_{\mathrm{de\;Rham}}
    \mathrm{Forms}".
\end{align*}
A construction of such functors (for any generalized cohomology theory $E$)
was given by Hopkins and Singer \cite[Def. 4.34]{hopkins2005}, and from the modern viewpoint
they can be described quite efficiently in a very general setting
via sheaves of spectra \cite{bunke2016differential}.

In order to understand and compute these abstractly defined refinements, it is
however important to have concrete models. Differential $K$-theory is
an especially
prominent example of this, since it appears in mathematical as well as physics
discussions, often in the form of a geometric model. In the case of $K$-theory,
the differential version is $\Z_2$-graded and the even and odd part were
developed independently. On the category of compact manifolds, a variety of
descriptions are available.
Simons and Sullivan \cite[\S 3]{simons_structured_2008} show that even differential
$K$-theory is defined by structured vector
bundles, i.e.\ vector bundles with connection with a suitable equivalence
relation. This picture was completed by Tradler, Wilson and Zeinalian
\cite[Thm. 5.7]{tradler_elementary_2012} by giving
a geometric description of odd differential $K$-theory via operator theory,
presented as
maps into the stable unitary group, where the addition is induced by a blocksum
operation. Later, via the Caloron correspondence, an interpretation of
their model via $\Omega$-bundles was developed in \cite[Thm.
3.17]{hekmati2015geometric}.

More recently, another approach has been implemented in
\cite[Thm. 4.25]{tradler_differential_2015}. The authors discuss the question of
representability of the $\hat K$-functor. As any cohomology theory, topological
$K$-theory is represented by homotopy classes of maps into the corresponding
spectrum, i.e. 
\begin{align*}
    K^0(M) \cong [M,B\mathrm U\times \Z], \qquad K^1(M) \cong [M,\mathrm U],
\end{align*}
where $\mathrm U$ is the stable unitary group, i.e. the union over all
$\mathrm U(n)$. For compact manifolds, this agrees with the usual description as the
Grothendieck group of the monoid of complex isomorphism classes of vector
bundles. For non-compact manifolds, we can take this as a definition (the vector
bundle definition would not yield a cohomology theory). Since only
the homotopy type of these spaces is relevant in this 
description, one can find good models for $B\mathrm U\times \Z$ and $\mathrm U$, which
carry the additional information needed to define a differential $K$-theory
class from a map into it. In the end, the authors describe
even and odd differential 
$K$-theory via smooth maps into explicit classifying spaces, equipped with
differential forms that represent the universal Chern character. These universal
forms are defined on filtrations of their spaces via compact smooth manifolds
(the usual finite-dimensional Grassmannians and unitary groups).
Therefore, again, this method relies heavily on the fact that a compact smooth
manifold will always map to a finite stage in the filtration. The problem with
working directly on the spaces $B\mathrm U\times \Z$
and $\mathrm U$ is of course their infinite-dimensional nature. As colimits of
finite-dimensional smooth manifolds, they are Fr\'echet manifolds, and as such,
it is harder to for example talk about differential forms on them. 

The new approach employed in this paper is the use of operator theory 
to perform certain norm completions and slightly enlarge these spaces in order
to improve their regularity. The result are well-behaved Banach manifolds, which
we then equip with natural differential forms in the classical sense. These
constructions are closely related to to the identification of
$B\mathrm U\times \Z$ with the space of Fredholm operators via a generalized
index map, as shown by Atiyah and Jänich.
While \cite{tradler_differential_2015} ultimately work with Chen spaces as models for the spaces
$B\mathrm U \times \Z$ and $\mathrm U$, our Banach manifolds allow us to do certain
calculations directly in the universal example,
without considering test manifolds. One immediate advantage of this approach is
that our model also works in the case of non-compact manifolds. If we restrict
to the compact case, we have explicit maps from the \cite{tradler_differential_2015} classifying spaces into
our completed versions, which induce isomorphisms in differential $K$-theory.

The addition map in $\hat K$ is implemented in both the odd and even case
via an explicit blocksum operation on the classifying spaces, which we denote by
$\boxplus$. Inverting an element corresponds to 
taking the operator adjoint in the odd case, and flipping the polarization on
the underlying polarized Hilbert space
in the even case. Since the \cite{tradler_differential_2015}-model operates always in a finite step in a filtration
of the classifying spaces by compact manifolds, the authors have to implement
certain finite-dimensional
shifts which depend on choices of a concrete representative of a cocycle. Our
model gets rid of the need for this shift by working directly in the
infinite-dimensional setting. Same as in their model, a map to our
geometrically enriched classifying spaces carries more data than just the
isomorphism class of a bundle, and is actually all that is needed to define a
differential $K$-theory class. In this sense, we take the classifying map
approach seriously and do not need the additional data of a
differential form that is present in other models. We also produce an obvious
cycle map, that assigns a differential $K$-theory class to a vector bundle with
connection. These are crucial differences to the abstract spectrum-based
construction given by Hopkins and Singer. Our main result is therefore
\begin{thm}
    Let $\grres$ be the restricted Grassmannian (Def.
    \ref{dfn:restrictedgrassmannian}) and $\mathrm U^1$ be the unitary group
    of operators which have a determinant (Def. \ref{dfn:u1}). On the category
    of possibly non-compact smooth manifolds, the abelian group valued functors 
    \begin{align*}
        \hat K^0(M) = \mathrm{Map}(M, \grres) / \CS\mathrm{-homotopy} +
        \mathrm{Stabilization}\\
        \hat K^{1}(M) = \mathrm{Map}(M, \mathrm U^1) / \CS\mathrm{-homotopy} +
        \mathrm{Stabilization}
    \end{align*}
    define differential $K$-theory.
\end{thm}
The equivalence relation is a geometrically refined version of homotopy (see
Def. \ref{dfn:CSequivalence})
that makes it possible to extract a differential form representative for the
Chern character out of an equivalence class\footnote{Note that this is
impossible for a homotopy class of maps.}, plus the
additional stability requirement that blocksumming with a constant map
to the basepoint does not change the equivalence class. It is a natural
question whether the stability relation is already contained in the
$\CS$-homotopy equivalence class, i.e. whether $f\sim_{\CS} f\boxplus \mathrm{const}_*$
for any representative of a $\hat K$-class. In the compact case, we get an
affirmative answer to this question and we prove
\begin{thm}
    \label{thm:mainthmintro}
    On the category of compact smooth manifolds, the abelian group valued functors 
    \begin{align*}
        \hat K^0(M) = \mathrm{Map}(M, \grres) / \CS\mathrm{-homotopy}\\
        \hat K^{1}(M) = \mathrm{Map}(M, \mathrm U^1) /\CS\mathrm{-homotopy}
    \end{align*}
    define differential $K$-theory.
\end{thm}
\setcounter{section}{1}
The idea that leads to this operator-theoretic approach can be described as follows: While 
$K$-theory is the study of stable vector bundles, it can also be interpreted as studying
Hilbert space bundles with a reduction of the structure
group to the stable general linear group $\mathrm{GL} \subset \mathrm{GL}(\h)$,
sitting in the
(contractible) full general linear group of $\h$. By Palais' tame approximation
theorem \cite[Thm. B]{PALAIS1965271} this
group is homotopy equivalent via its natural inclusion to the group of operators
which have a determinant, denoted by $\mathrm{GL}^1$. Therefore we might as well
study the space $B\mathrm{GL}^1$. There happens to be a model of the universal
smooth principal $\mathrm{GL}^1$-fiber bundle, which has appeared in the study of loop
groups \cite[Sec. 7.5]{pressley_loop_1988} and also in applications in physics
in the form of fermionic
second quantization (for a mathematical treatment see
\cite[Sec. V.2]{wurzbacher2001fermionic}). This bundle carries a connection,
which gives rise to a universal
Chern character differential form via the usual Chern-Weil formula. The degree
$2$-part of this form is known in
the physics literature as the Schwinger cocycle, where the
discussion usually focusses on line bundles. We prove
that we can get representatives also for the higher dimensional parts of the
Chern character (as has already been observed by \cite[Thm. 3.9]{freed_index_1988}), and
along the way, we review some constructions in the world of
restricted unitary groups, Grassmannians and Stiefel manifolds, which we could
not find a good reference for. 

Our proof proceeds in the following way. In Section \ref{sec:universal} and
\ref{sec:restricted}, we review the construction of the restricted Grassmannian
and the unitary group of operators which have a determinant, which will give the
even resp. odd model for differential $K$-theory. The universal Chern class in
the odd case is induced by the Maurer-Carten form of $\mathrm U^1$. In the even
case, we review the construction of a certain universal bundles over $\grres$, the
curvature of which gives rise to invariant representatives of the Chern
character via Chern--Weil theory.

In Section \ref{sec:CSformsandblocksum}, we equip these spaces with an $H$-space
structure. The key difference to the purely homotopy theoretical
approach is that we have to choose these structures in such a way that they are
compatible with the Chern and Chern--Simons forms. For example, even though it
induces addition in $K$-theory, operator
multiplication on the unitary group will not work as an addition in $\hat K^{1}$, since it
will not make the Chern character map into a monoid morphism on the level of
differential forms.

Section \ref{sec:smoothper} reviews geometric versions of the usual periodicity
maps in the $K$-theory spectrum. The even to odd part is given by the holonomy
map in the universal fibration, while the odd to even part is a certain
multiplication operator map considered already by Pressley and Segal in their
study of loop group representations \cite[Sec. 6.3]{pressley_loop_1988}. It is
interesting, though not a key fact for us, that this map can be used to
implement Bott periodicity as a smooth homomorphism of infinite-dimensional
Lie groups. We also prove that the geometric spaces we use combine to an
$\Omega$-spectrum representing $K$-theory, where the addition is implemented by
our blocksum (Prop. \ref{prop:spectrum}).

In Section \ref{sec:differential}, we put together all the ingredients from the
previous sections in order to prove that the
previously discussed blocksum and inversion operations equip the Chern--Simons
equivalence classes of maps into the classifying spaces with an abelian group
structure. This is achieved by finding
explicit homotopies directly on our classifying spaces, which need to have
vanishing Chern--Simons forms. The discussion here is simplified considerably by
the simple cohomological structure of the relevant spaces and the availability
of a de Rham theorem for the Banach manifolds in question\footnote{Since not all
Banach manifolds admit a smooth partition of unity, this is not immediately
obvious.}.

Having the abelian group structure on $\hat K^0$ and $\hat K^1$, what is left to
do in Section \ref{sec:natural} is to give the remaining structure maps for a
differential extension and check the corresponding axioms. Here, the
periodicity maps constructed in Section \ref{sec:smoothper} play a key role.

In Section \ref{sec:comparison} we make a comparison to the
\cite{tradler_differential_2015}-model
of differential $K$-theory and find isomorphisms induced by explicit maps
on the classifying spaces. From this isomorphism we learn that the
stability condition in our main theorem can actually be dropped if we 
restrict to compact manifolds. We close by discussing some examples of
differential $K$-theory classes in Section $\ref{sec:examples}$.

\textbf{Acknowledgements.} I thank the City University of New York,
especially Scott O. Wilson, for their hospitality during my research
stay, where parts of this research were conducted. I also thank Mahmoud
Zeinalian and Matthew Cushman for many useful conversations, and Mauricio
Bustamante and Markus
Upmeier for valuable suggestions on an earlier draft of this paper. This work is
part of the author's PhD thesis at the University of Augsburg under the supervision
of Bernhard Hanke.

\section{Universal representatives for the Chern character}
\label{sec:universal}

Central to this work are the constructions of explicit smooth models
for the classifying spaces of even and odd $K$-theory. Recall that the complex
$K$-theory spectrum is two-periodic and consists of the
spaces $B\mathrm U\times \Z$ in the even degrees and $\mathrm U$ in the odd
degrees, where $\mathrm U$ is the stable unitary group, i.e. the colimit along
the inclusions $\mathrm U(n)\hookrightarrow \mathrm U(n+1)$. In order to
build a differential extension of $K$-theory, we define smooth models for both
of these spaces which carry natural invariant differential forms that represent the
universal Chern character.

For the odd case, recall that on $\mathrm U(n)$, we have the Maurer-Cartan form
$\omega_n$. It is well known that the real cohomology of $\mathrm U(n)$ is generated by
the cohomology classes represented by the invariant differential forms
\begin{align}
    \label{eq:MCform}
    \left(\frac{i}{2\pi}\right)^{k} \frac{(-1)^{k-1}(k-1)!}{(2k-1)!} \tr \left( 
    \omega_n^{2k-1}\right)\in \Omega^{2k-1}(\mathrm U(n)). 
\end{align}
The normalizations we have chosen here are in order to make this agree with the
transgression of the Chern character in the universal fibration (see Section
\ref{sec:smoothper}). We can stabilize using the usual inclusion $\mathrm U(n)
\hookrightarrow \mathrm U(n+1)$, but when one goes
to limit, one has to deal with the intricacies of infinite-dimensional
manifolds. Our preferred way of dealing with this is to work in the setting of
Banach manifolds. The problem is that the Lie algebra of the stable unitary group
$\mathrm U$ is supposed to consist of skew-adjoint finite rank operators of
arbitrary dimension. Since this
is not a closed subspace of the bounded operators, there are some complications
if we want to consider $\mathrm U$ as a smooth manifold. A simple fix is
to instead go one step further and complete with respect to the trace norm
\begin{align*}
    || X ||_{L^1} = \tr |X| = \tr \sqrt{X^* X}.
\end{align*}
This leads to the ideal $L^1$ of trace-class operators, and further to the
Banach-Lie group $\mathrm U^1$, which we now define. 
\begin{dfn}
    \label{dfn:u1}
    Let $\h$ be a complex
    separable infinite-dimensional Hilbert space. Then $\mathrm U^1$ is the subgroup of the
    unitaries of $\h$ given by
    \begin{align*}
        \mathrm U^1 = \left\{ P\in \mathrm U(\h) \mid P-1\in L^1 \right\},
    \end{align*}
    with topology induced by the inclusion
    \begin{align*}
        \mathrm U^1&\hookrightarrow L^1\\ P &\mapsto P-1.
    \end{align*}
\end{dfn}
Palais \cite[Thm. B]{PALAIS1965271} showed that the inclusion of the stable unitary group $\mathrm
U\hookrightarrow\mathrm U^1$ is a homotopy
equivalence, but $\mathrm U^1$ has better regularity, as it is actually a Banach-Lie
group, locally modelled on the Banach space $L^1$. It is well known that its
cohomology is generated entirely by traces of odd powers of the Maurer-Cartan
form $\omega$, analogous to formula (\ref{eq:MCform}). It is therefore sensible
to make the following definition.
\begin{dfn}
    \label{dfn:oddchern}
    The universal odd Chern character form $\ch_{\mathrm{odd}}\in
    \Omega^{\text{odd}}(\mathrm U^1)$ is 
    \begin{align*}
        \ch_{\mathrm{odd}} = \sum_{k\geq 1} \ch_{2k-1} = \sum_{k\geq 1}
        \left(\frac{i}{2\pi}\right)^{k}\frac{(-1)^{k-1}(k-1)!}{(2k-1)!}
        \tr \left(\omega^{2k-1}\right).
    \end{align*}
\end{dfn}
In order to find a good model for the even case, we recall the construction of
universal connections. We will first review the situation for the
finite-dimensional Grassmannians, and then spend the next chapter to generalize
to the infinite-dimensional setting. As one would hope, these universal
connections will yield well suited differential form representatives for the
universal Chern character on our Grassmannian model of $B \mathrm U \times
\Z$. 

The Stiefel bundle over the Grassmannian manifold is given by
\begin{align*}
    \mathrm{St}_{k,N} = \mathrm U(N)/I_k \times \mathrm U(N-k) \to \mathrm U(N)
    / \mathrm U(k)\times \mathrm U(N-k) = \mathrm{Gr}_{k,N}. 
\end{align*}
There is a map $S\colon \mathrm{St}_{k,N} \to M_{N\times k}$ which assigns to an element
on the left a matrix $A\in M_{N\times k}$ which
satisfies $A^* A = I_k$. The entries of $A$ are just given by the first
$k$ columns of a representative of our left coset. Denote by $S^*$ the map
$S$ followed by taking the adjoint matrix, and denote by $\dd S$ the
differential of $S$, which is an $M_{N\times k}$-valued differential form. Then,
there is a Lie algebra valued $1$-form given by $S^*\dd S$, and one can show
that it takes values in the skew adjoint matrices and furthermore that it defines 
a connection for the given
principal bundle. Narasimhan and Ramanan \cite[Thm. 1]{narasimhan_existence_1961}
observed that
the family of connections given by this construction for varying $k$ and
$N$ have a universal property, meaning that every smooth principal bundle for a
unitary group with a given connection comes from pulling back such a bundle
and its respective connection by a smooth classifying map.

By Chern--Weil theory, one can define representatives for the Chern character
by chosing a connection and considering traces of powers of its curvature. The curvature of
$\omega=S^*\dd S$ can be calculated explicitely as follows. An element in the
tangent space at $\begin{pmatrix}
    I_k \\ 0
\end{pmatrix}$ of $\mathrm{St}_{k,N}$ is given by an $n\times n$ skew-hermitian block
matrix $\begin{pmatrix}
        P & -Q^* \\ Q & 0 
    \end{pmatrix}$
where $P$ is a skew hermitian $k\times k$ matrix and $Q$ is an arbitrary
$(N\times (N-k))$ - matrix. The horizontal subspace is given by the kernel of
$\omega$, which corresponds to matrices which have $P=0$. Recall that the
curvature according to \cite[Thm. 5.2]{kobayashi1996foundations}
is defined to be the covariant derivative of the connection, so we
have $\Omega = \dd \omega \circ h$, where $h$ is the horizontal projection. We
calculate
\begin{align}
    \label{eq:curvature}
    \Omega\left(\begin{pmatrix}
        P_1 & -Q_1^* \\ Q_1 & 0 
    \end{pmatrix}, \begin{pmatrix}
        P_2 & -Q_2^* \\ Q_2 & 0 
    \end{pmatrix}\right) &= \dd \omega \left( \begin{pmatrix}
        0 & -Q_1^* \\ Q_1 & 0 
    \end{pmatrix}, \begin{pmatrix}
        0 & -Q_2^* \\ Q_2 & 0 
    \end{pmatrix} \right) \nonumber\\&=-\omega \left[ \begin{pmatrix}
        0 & -Q_1^* \\ Q_1 & 0 
    \end{pmatrix}, \begin{pmatrix}
        0 & -Q_2^* \\ Q_2 & 0 
    \end{pmatrix} \right] \nonumber \\&= Q_1^* Q_2 -Q_2^* Q_1.
\end{align}
Invariance under the transitive left $\mathrm U(N)$-action allows us to
extend this to at any point in $\mathrm{St}_{k,N}$. The usual Chern--Weil theory
then gives explicit differential forms on the Grassmannian after we take traces.

As in the odd case, these invariant forms stabilize under the inclusions
$\mathrm{Gr}_{k,N}
\hookrightarrow \mathrm{Gr}_{k,N+1}$, but again, when we want to work with a 
universal space, problems arise. The direct limit of the Grassmannians is not a
Banach manifold, and so one needs more delicate tools to talk about connections
and even differential forms on them. There is no obvious construction of a universal
invariant connection for $\mathrm U$-bundles in the
stable case, and some of the problems that arise are discussed in
in \cite[Prop. 2.3]{freed_index_1988}. However there still exists an 
analog to the finite-dimensional construction in the category of Banach
manifolds, which we will review in the next section.

\section{The restricted Stiefel manifold and Grassmannian}
\label{sec:restricted}

In the infinite-dimensional
setting, for a Hilbert space $\h$, the unitary group $\mathrm U(\h)$ becomes
contractible, so one usually restricts to appropriate subgroups in order to
generate non-trivial topology. Assume that our Hilbert space $\h$ (complex,
separable, infinite-dimensional) comes with a $\Z$-graded orthonormal basis $\left\{ e_i
\right\}_{i\in \Z}$, thereby defining a grading (also sometimes called
polarization) into two infinite-dimensional, complementary subspaces 
\begin{align*}
    \h\cong \h_+\oplus \h_- = \mathrm{span}\left\{ e_i \mid i\geq 0
    \right\}\oplus \mathrm{span}\left\{ e_i \mid i < 0 \right\}.
\end{align*}
The grading can also be seen as given by the involution $\varepsilon = \begin{pmatrix}
    1&0\\0&-1
\end{pmatrix}$. We define the Banach algebra of bounded operators
\begin{align*}
    \lglres = \left\{ \begin{pmatrix}
        X_{++} &X_{-+} \\X_{+-}&X_{--}
\end{pmatrix}\in \mathfrak{gl}(\h_+\oplus \h_-) \mid X_{-+},X_{+-}\in L^2 \right\}
\end{align*}
with norm
\begin{align*}
    ||\begin{pmatrix} X_{++} &X_{-+} \\X_{+-}&X_{--} \end{pmatrix} || =
    ||X_{++}|| + ||X_{--}|| + ||X_{-+}||_{L^2}
    + ||X_{+-}||_{L^2}.
\end{align*}
Recall that $L^2$ denotes the ideal of Hilbert-Schmidt operators, i.e. operators
that meet the summability condition $\mathrm{tr}\,X^*X  < \infty$. One could
equivalently define $\lglres$ to be the subalgebra of bounded operators that
commute with $\varepsilon$ up to a Hilbert-Schmidt operator. The group of
units in this Banach algebra is the restricted general linear group
$\glres$ of \cite[Sec. 6.2]{pressley_loop_1988}. It is easy to see
that for $\begin{pmatrix} X_{++} &X_{-+} \\X_{+-}&X_{--} \end{pmatrix}\in \glres$, the operators
$X_{++}$ and $X_{--}$ have to be Fredholm operators, since they are invertible up to
compacts. Then, one can show that the projection
\begin{align}
    \psi\colon\glres \to \text{Fred}\nonumber\\
    \begin{pmatrix} X_{++} &X_{-+} \\X_{+-}&X_{--} \end{pmatrix} \mapsto X_{++}
    \label{eq:htpyeq}
\end{align}
is a homotopy equivalence \cite[Cor. 3.1]{wurzbacher_elementary_2006}. By polar
decomposition, the restricted unitary group $\ures= \glres \cap \mathrm U(\h)$
is homotopy equivalent to
$\glres$, and since by the Atiyah-Jänich theorem, $\mathrm{Fred}\sim
B\mathrm U\times \Z$, this makes $\ures $ into a suitable Banach-Lie
group model for $B\mathrm U\times \Z$. We will now consider the associated
Grassmanian to this situation.
\begin{dfn}
    \label{dfn:restrictedgrassmannian}
    The restricted Grassmannian $\grres$ is the set of
    all closed subspaces $W\subset \h$ such that the orthogonal projection
    $\pi_+\colon W\to \h_+$ is a Fredholm operator and $\pi_-\colon W\to \h_-$ is a
    Hilbert-Schmidt operator.
\end{dfn}
Loosely speaking, we only consider subspaces here which are
comparable in size with $\h_+$, in the sense of a perturbation by a
Hilbert-Schmidt operator. As in the finite-dimensional case,
there are many equivalent descriptions of the Grassmannian. 
\begin{prop}
    \label{prop:descriptionsofgrres}
    A point in
    $\grres$ can be thought of as
    \begin{enumerate}[(i)]
        \item A subspace $W\subset \h$ such that $\restr{\pi_+}{W}\in \fred$ and
            $\restr{\pi_-}{W}\in L^2$.
        \item A self-adjoint projection operator $\pi$ on $\h$ such that the
            commutator $[\pi,\varepsilon]\in L^2$.
        \item A self-adjoint involution $F$ on $\h$ such that $F-\varepsilon\in
            L^1$.
        \item An equivalence class $[X]\in \ures / \mathrm{U}(\h_+) \times
            \mathrm{U}(\h_-)$.
    \end{enumerate}
    \begin{proof}
        In \cite[Prop. 7.1.3]{pressley_loop_1988}, it is proved that $\ures$ acts
        transitively on $\grres$ with stabilizer $\mathrm{U}(\h_+) \times
            \mathrm{U}(\h_-)$ and thus we get the description (iv) as a
        homogenous space. To a representative $X\in \ures$, we associate the subspace
        $W=X(\h_+)$ to get back to (i). Furthermore, we can construct the self
        adjoint involution $F= X\varepsilon X^*$, which is $+\id$ on the subspace $W$ and
        $-\id$ on $W^{\perp}$, to get to (iii). Lastly, we can consider the
        projection operator $(F+\id)/2$, which gives (ii).
    \end{proof}
\end{prop}
It is often convenient to have multiple descriptions of
$\grres$. Note that using $(\mathrm{iv})$, we can endow $\grres$ with the structure of a
Hilbert manifold modelled on 
\begin{align*}
    T_{1} \grres \cong \lures / \mathfrak{u}(\h_+)\times \mathfrak{u}(\h_-)
    \cong L^2(\h_+,\h_-).
\end{align*}
By Kuiper's theorem, $\mathrm U(\h_{\pm})$ is contractible. Since the projection 
map $\ures\to \grres$ defines a locally trivial principal bundle, it is
therefore actually a homotopy equivalence, in sharp
contrast to the finite-dimensional case (for details, see
\cite[Lemma 2.1]{wurzbacher_elementary_2006}). It follows that the restricted Grassmannian 
has infinitely many
diffeomorphic path components, indexed by $\Z$, which can be recovered from a
given subspace $W$ by its virtual dimension
\begin{align*}
    \mathrm{virt.dim}(W) = \mathrm{dim}(\mathrm{ker}(\pi_+\colon W\to \h_+))-
    \mathrm{dim}(\mathrm{coker}(\pi_+\colon W\to \h_+)).
\end{align*}
If $W = X(\h_+)$ for $X\in \ures$, then $\mathrm{virt.dim}(W) =
\mathrm{ind}(X_{++})$.
As in the finite-dimensional case, there is a corresponding Stiefel manifold. 
\begin{dfn}
    \label{dfn:stiefelmfd}
    The restricted Stiefel manifold is the set of operators
    \begin{align*}
        \mathrm{St}_{\mathrm{res}} = \left\{ w = \begin{pmatrix}
            w_+ \\ w_-
        \end{pmatrix}\in \mathfrak{gl}(\h_+,\h) \mid w \;\mathrm{injective}, w_+ - 1 \in L^1, w_-\in
    L^2\right\},
    \end{align*}
    endowed with the topology and smooth structure coming from the inclusion as an
    open subset into the Banach space $L^1\times L^2$.
\end{dfn}
An element in $\stres$ is sometimes called an admissible base. We have the
following immediate observation.
\begin{prop}
    \label{prop:contractible}
    The restricted Stiefel manifold is contractible.
    \begin{proof}
        Consider the filtration of $\h_+$ by the finite-dimensional subspaces
        $V_N$ spanned by the first $N$ basis vectors. By Palais' tame
        approximation theorem \cite[Thm. A]{PALAIS1965271}, $\stres$ is
        homotopy equivalent to the inductive limit of the subspaces
        \begin{align*}
            \stres^N = \left\{ w\in \stres \mid \; w \pi_{V_n^\perp} = 0
        \right\}.
        \end{align*}
        But those are just the usual $n$-dimensional Stiefel manifolds that
        serve as the total space of the universal $\mathrm{GL}(n)$-bundle and as
        such are contractible.
    \end{proof}
\end{prop}
We have set up a situation very similar to the finite-dimensional one, where
one has a principal $\mathrm U(k)$-bundle $\mathrm{St}_{k,N}\to
\mathrm{Gr}_{k,N}$. The correct structure group in our case is the group of 
invertible operators which have a determinant
$\mathrm{GL}^1 = \left\{ P\in \mathrm{GL}(\h_+) \mid P-1\in L^1 \right\}$. This
group acts on $\stres$ on the right via the usual change of basis $(w,Q)\mapsto
w Q$. With this action, we have the following
\begin{prop}
    \label{prop:principal}
    The map $q\colon\stres\to\grres^0$ defines a smooth principal
    $\mathrm{GL}^1$-bundle over the path component of the basepoint $\h_+$ on
    the restricted Grassmannian.
    \begin{proof}
        The action is smooth since it is just multiplication of operators, and
        it is also clear that it is free. For fiberwise transitivity, we need to
        check that two admissible bases
        for the same subspace are related by right multiplication with elements
        in $\mathrm{GL}^1$. Let
        $w,w'$ be two admissible bases for $W$. Then
        $w' = w Q$, where $Q = w^{-1} w'\in \mathrm{GL}(\h_{+})$ and we need to
        show that $Q\in \mathrm{GL}^1$. We calculate
        \begin{align*}
            1+ L^1 = \pi_+ w' = \pi_+ w Q = Q + L^1.
        \end{align*}
        The only thing left to show is local triviality. As in the finite
        dimensional case, there exist graph coordinates for
        the restricted Grassmannian (cf. \cite[Ch. 7]{pressley_loop_1988}).
        Those are given by 
        \begin{gather*}
            L^2(W,W^{\perp})\to U\subset \grres\\T \mapsto \Gamma_T = \left\{
            (v,Tv) \mid v\in W \right\}.
        \end{gather*}
        Choose an $X\in \ures$ such that $W = X(\h_+)$. Then
        we define a local section by setting $s(T) = \restr{X}{\h_+} + T
        \restr{X}{\h_+}\in \stres$.
    \end{proof}
\end{prop}
\begin{rem}
    \label{rem:nohomogenousconnection}
    One would like to reduce the structure group of this bundle to the
    unitary group $\mathrm U^1$. Interestingly, this is actually not possible,
    since it would determine a homogenous connection which would ultimately
    imply that the bundle is trivial. This is discussed after Prop. 3.15 in
    \cite{freed_index_1988}. 
\end{rem}
\begin{cor}
    \label{cor:universalbundle}
    The smooth fiber bundle of Banach manifolds
    \begin{align*}
        \mathrm{GL}^1 \to \stres \to \grres^0
    \end{align*}
    is a model for the universal $\mathrm{GL}^1$-fibration.
\end{cor}
We will now construct a connection form for this
principal bundle that is supposed to represent the limit of the finite
dimensional connections on the bundles $\mathrm{St}_{k,N}\to \mathrm{Gr}_{k,N}$.
It will in particular generate representatives for the Chern
character which are compatible with the finite-dimensional versions. Consider
the coordinate map
\begin{align*}
    w\colon \stres &\to L^1\times L^2\\
    \begin{pmatrix}
        w_+ \\w_-
    \end{pmatrix}&\mapsto \begin{pmatrix}
        w_+-1\\ w_-
    \end{pmatrix},
\end{align*}
and consider its differential $\dd w$ as an operator-valued differential form on
$\stres$. Furthermore, we can associate to $w\in\stres$ the projection operator
$\pi_W\in \grres^0$ onto $W =w(\h_+)$, which gives another operator-valued differential
form $\dd \pi_W$ on $\stres$.
\begin{prop}
    \label{prop:univconnection}
    The assignment $\Theta = w^{-1} \pi_{W} \dd w$ defines a principal connection
    on $\mathrm{GL}^1$-bundle $\stres \to \grres^0$. The curvature of
    $\Theta$ is given by the expression 
    \begin{align*}
        \Omega = \dd \Theta + \frac{1}{2} [\Theta, \Theta] = w^{-1} \pi_{W} \dd
        \pi_W \dd \pi_W w.
    \end{align*}
    \begin{proof}
        We first check that $\Theta$ is $L^1$-valued. We can write $\Theta =
        w^{-1} \pi_W (\pi_+ + \pi_-) \dd w$, and since $\pi_+\dd w$ is trace
        class, it remains to show that the second summand is also trace class. A
        simple calculation shows that $\pi_W\in \grres$ is equivalent to
        $\pi_+-\pi_W\in L^2$. Therefore, using that $\pi_- \dd w = \dd (\pi_- w)
        \in L^2$, we have
        \begin{align*}
            w^{-1}\pi_W \pi_- \dd w = w^{-1} (\pi_+ + L^2) \pi_- \dd w = 0 + L^2
            \cdot L^2 =  L^1.
        \end{align*}
        We now check the defining properties of a connection form. On the
        fundamental vector fields for $X\in L^1$ of the form
        $\widetilde{X}_w = \restr{\frac{\dd}{\dd t}}{t=0} w \exp(t X)$, we
        clearly have $\Theta(\widetilde{X})=X$. On the other hand, we have
        \begin{align*}
            (R_Q^*{\Theta})_w = (wQ)^{-1} \pi_W (\dd w) Q = \mathrm{Ad}_{Q^{-1}}
            \Theta_w,
        \end{align*}
        finishing the proof that $\Theta$ is a connection form.
        
        For the calculation of the curvature, we will need the identities 
        \begin{gather*}
            \dd w = \dd (\pi_W w) = \dd \pi_W w + \pi_W \dd w \\ \dd \pi_W = \dd
            (w w^{-1} \pi_W) = \dd w w^{-1} \pi_W + w \dd (w^{-1}\pi_W).
        \end{gather*}
        From the second identity, it follows that 
        \begin{align*}
            \dd (w^{-1}\pi_W) = w^{-1}\pi_W \dd \pi_W - w^{-1}\pi_W\dd w w^{-1}
            \pi_W.
        \end{align*}
        We calculate
        \begin{align*}
            \dd \Theta = \dd (w^{-1} \pi_W \dd w) &= \dd (w^{-1} \pi_W) \dd w
            \\&= (w^{-1}\pi_W \dd \pi_W - w^{-1}\pi_W\dd w w^{-1} \pi_W) \dd w
            \\&= w^{-1} \pi_W \dd \pi_W (\dd \pi_W w + \pi_W \dd w) -
            w^{-1}\pi_W \dd w w^{-1} \pi_W \dd w  \\ &= w^{-1} \pi_W \dd \pi_W
            \dd \pi_W w - \frac{1}{2}[\Theta,\Theta],
        \end{align*}
        since $\pi_W \dd\pi_W \pi_W=0$.
    \end{proof}
\end{prop}
Since the curvature form is trace class valued, the usual arguments from
Chern--Weil theory go through and give representatives for the Chern character of
the universal $\mathrm{GL}^1$ bundle over $\grres$ (cf. \cite[Thm.
1.13]{freed_index_1988}). One difference to the
bundles over the finite-dimensional Grassmannians is that our form $\Theta$ is
not left-invariant under the action of $\mathscr{E}$, since the left action of
$(1,Q)\in \mathrm{GL}^1$ is the same as acting on the right by $Q^{-1}$ - an
operation that has to be equivariant with respect to the adjoint representation,
as checked in the proof of Prop. \ref{prop:univconnection}. However, we still
have that after taking traces, the forms $\tr\, \Omega^k$ make sense as 
invariant differential forms on $\grres^0$, which invariantly extend to the
other diffeomorphic components of $\grres$. We make the following definition
\begin{dfn}
    \label{dfn:evenchernform}
    The univeral even Chern character form
    $\ch_{\mathrm{even}}\in\Omega^{\mathrm{even}}(\grres)$ is 
    \begin{align*}
        \ch_{\mathrm{even}} = \sum_{k\geq 0} \ch_{2k} = \ch_0 + \sum_{k\geq 1}
        \left( \frac{i}{2 \pi} \right)^k
        \frac{1}{k!} \tr \left(\Omega^{k}\right),
    \end{align*}
    where $\Omega = \pi_W \dd \pi_W \dd \pi_W$ is a trace class operator valued
    form. Here, $\ch_0\colon \grres \to \Z$ is the map that assigns to $W$ its
    virtual dimension. 
\end{dfn}
The positive degree forms are actually invariant: Since the
action of $\ures$ is by conjugation of both $\pi_W$ and $\dd \pi_W$ by a unitary, it
leaves the trace invariant. Thus, it is useful to 
explicitly work out what happens at the tangent space of $\h_+$. Recall that
\begin{align*}
    T_1 \grres \cong \lures / \mathfrak{u}(\h_+) \times \mathfrak{u}(\h_-)\cong
    \left\{ \begin{pmatrix}
        0&-c^* \\ c&0
\end{pmatrix} \mid c\in L^2(\h_+,\h_-)\right\}.
\end{align*}
Set $w_0 = \begin{pmatrix}
    1 \\ 0
\end{pmatrix}\in \stres$. For $w=X w_0\in \stres$, we have that $\pi_W = \pi_{X(\h_+)} =
X\pi_{+}X^*$ and therefore $(\dd \pi_W)_{\pi_+} = [-, \pi_{+}]$, where the bracket
indicates the commutator. Therefore, evaluation of $\Omega = \pi_W \dd \pi_W \dd
\pi_W$ at the point $\pi_{+}$ yields
\begin{gather}
    \label{eq:univcurvature}
    \Omega_{\pi_{+}} (\begin{pmatrix}
        0&-c_{1}^* \\ c_{1}&0
\end{pmatrix},\begin{pmatrix}
        0&-c_{2}^* \\ c_{2}&0
    \end{pmatrix}) \\=\begin{pmatrix}
        1&0 \\0 &0
    \end{pmatrix}\bigg(\left[ \begin{pmatrix}
            0&-c_{1}^* \\ c_{1}&0\end{pmatrix}, \begin{pmatrix}
                1&0\\0&0
            \end{pmatrix}
\right]\left[ \begin{pmatrix}
            0&-c_{2}^* \\ c_{2}&0\end{pmatrix}, \begin{pmatrix}
                1&0\\0&0
        \end{pmatrix} \right] -\nonumber \\\left[ \begin{pmatrix}
            0&-c_{2}^* \\ c_{2}&0\end{pmatrix}, \begin{pmatrix}
                1&0\\0&0
            \end{pmatrix}
\right]\left[ \begin{pmatrix}
            0&-c_{1}^* \\ c_{1}&0\end{pmatrix}, \begin{pmatrix}
                1&0\\0&0
        \end{pmatrix} \right]\bigg) \nonumber\\=\begin{pmatrix}
            1&0\\0&0
        \end{pmatrix} \left(\begin{pmatrix}
            0& c_1^* \\c_1&0
        \end{pmatrix} \begin{pmatrix}
            0&c_2^*\\c_2&0
        \end{pmatrix}- \begin{pmatrix}
            0& c_2^* \\c_2&0
        \end{pmatrix} \begin{pmatrix}
            0&c_1^*\\c_1&0
        \end{pmatrix}  \right) \nonumber\\= \begin{pmatrix}
            1&0\\0&0
        \end{pmatrix} \begin{pmatrix}
            c_1^* c_2 - c_2^* c_1 & 0 \\ 0 & c_1 c_2^* - c_2 c_1^*
        \end{pmatrix} = \begin{pmatrix}
            c_1^* c_2 - c_2^* c_1 & 0 \\ 0 & 0
        \end{pmatrix},\nonumber
\end{gather}
and we recover the familiar formula from the finite-dimensional case
(\ref{eq:curvature}). 

There are natural smooth inclusions of the finite-dimensional
Grassmannians into the restricted Grassmannian, given as follows: Pick a
$\Z$-graded orthonormal basis $\left\{ e_i \right\}$ for $\h$, where
$\h_+\cong\mathrm{span}\left\{ e_i \mid i \geq 0\right\} $. Considering for
$N\in\Z$ the subspaces 
\begin{align*}
    \h_N = \mathrm{span} \left\{ e_i  \mid i \geq N\right\},
\end{align*}
one sees that the subsets
\begin{align*}
    \mathrm{Gr}_{\mathrm{res},N} = \left\{ W\in \grres \mid \h_{N}\subset W\subset \h_{-N} \right\}
\end{align*}
are isomorphic to the full finite-dimensional Grassmannians $\mathrm{Gr}(\C^{2N}) =
\coprod_{k\leq 2N} \mathrm{Gr}_{k,2N}$ by mapping $W$ to $W/\h_{N}
\subset\h_{-N}/\h_{N}\cong \C^{2N}$. The inclusion of
$\mathrm{Gr}_{\mathrm{res},N}$ into $\mathrm{Gr}_{\mathrm{res},N+1}$ corresponds
to sending $V\in \mathrm{Gr}(\C^{2N})$ to $\left\{ 0 \right\} \oplus V\oplus
\C\in \mathrm{Gr}(\C^{2(N+1)})$. The union of these finite-dimensional
Grassmannians, denoted by $\mathrm{Gr}_{\mathrm{res},\infty}$, is dense in
$\grres$, and the intersection $\mathrm{Gr}_{\mathrm{res},N}\cap \grres^k$ is
diffeomorphic to $\mathrm{Gr}_{N+k,2N}$ (cf. \cite[Prop.
III.5]{wurzbacher2001fermionic}). All in all, we have inclusion maps
\begin{align*}
    i\colon \mathrm{Gr}_{k,2N} = \mathrm{Gr}_{N+(k-N), 2N} &\to
    \mathrm{Gr}_{\mathrm{res},\infty}\subset \grres \\W &\mapsto
    W\oplus \h_N,
\end{align*}
which are easily seen to be compatible with the chosen Chern character
differential forms in the following sense: 
\begin{prop}
    \label{prop:compatbilityChern}
    Under the natural inclusion $i\colon \mathrm{Gr}_{k,2N} \hookrightarrow
    \grres$, the universal Chern character form $\ch_{\mathrm{even}}$ pulls back to
    the corresponding forms on the finite-dimensional Grassmannian, which are
    given by the Chern--Weil forms of the universal connection (see
    (\ref{eq:curvature})).
    \begin{proof}
        On the level of projections, with the
        above mentioned identification of $\C^{2N}$ with a subset of $\h$, we see
        that $\pi_W$ gets mapped by $i$ to $\pi_W + \pi_N$, where $\pi_N$ is the
        projection to $\h_N$. But that means that the pullback of the form can be
        written as
        \begin{align*}
            (\pi_W+ \pi_N) \dd (\pi_W+ \pi_N) \dd (\pi_W+ \pi_N) &= (\pi_W + \pi_N)
            \dd \pi_W \dd \pi_W  \\&= \pi_W \dd \pi_W \dd \pi_W  - \dd \pi_N \pi_W \dd
            \pi_ W \\&= \pi_W \dd \pi_W \dd \pi_W,
        \end{align*}
        where we used that $\pi_N \pi_W = \pi_W \pi_N = 0$. We can now do the
        calculation (\ref{eq:univcurvature}) again, but this time on the
        finite-dimensional Grassmannian, and see that (using additivity of the
        trace) the resulting Chern character forms $\tr\, i^*(\Omega^k)$ are precicely
        the ones coming from the Narasimhan--Ramanan curvature form.
    \end{proof}
\end{prop}
\begin{rem}
    \label{rem:cyclemap}
    We can use this calculation to cook up a ``cycle map'': Given a connected manifold
    $M$ and a class in $\hat K^0(M)$ represented by a formal difference
    $[V,\nabla_V] - [W,\nabla_W]$ of smooth hermitean vector bundles with
    compatible connections of dimension $k$ and $k'$, we can 
    use the Narasimhan--Ramanan theorem to get classifying maps $f_{V}\colon M\to
    \mathrm{Gr_{k,2N}}$, $f_{W}\colon M\to \mathrm{Gr_{k',2N}}$.
    Employing our above defined inclusions, we may as well assume
    that the target of these maps is actually $\grres$. Then, using the flip and
    blocksum map defined in Section \ref{sec:CSformsandblocksum}, we get a smooth map to
    the restricted Grassmannian, given by $f_V \boxplus \flip(f_W)$, which
    represents the differential $K$-theory class in our model.
    Note that $\ch_0(f_V \boxplus \flip(f_W)) = (k-N) - (l-N) = k-l =
    \ch_0(V) - \ch_0(W)$, which justifies our definition of the degree zero part $\ch_0$
    of the Chern character.
\end{rem}

\section{Chern--Simons forms, the blocksum and the inversion operation}
\label{sec:CSformsandblocksum}

We begin this chapter by discussing the transgressions of the Chern character in
the path loop fibration. The resulting Chern--Simons forms have first appeared in
\cite[Sec. 3]{chern1974characteristic} and they were one of the key ideas that lead to
the development of differential cohomology theories.

Let us consider the universal situation of the smooth path-loop fibration 
over $\mathrm U^1$ and $\grres$. There
are some subtleties when one wants to consider path and loop spaces as smooth
manifolds, but all we need is to have well-defined pullbacks to
finite-dimensional manifolds. This situation can be made precise by Chen's notion
of diffeological spaces \cite[Def. 1.2.1]{chen1977}. However, the 
identities that we want are provable via topological arguments, so this
viewpoint is not too important for the present paper, and one might as well
interpret the next paragraph as an informal motivation for the second part of Def.
\ref{dfn:chernmaps}.

By pulling back along the evaluation maps $P \grres \times I \to \grres$ and $P
\mathrm U^1 \times I \to \mathrm U^1$ and then fiber integrating, we
arrive at the universal Chern--Simons forms 
\begin{align*}
    \cs_{\mathrm{odd}} &= \int_{I} \mathrm{ev}_t^* (\ch_{\mathrm{even}}) \in
    \Omega^{\mathrm{odd}}(P U^1) \\\cs_{\mathrm{even}} &= \int_{I}
    \mathrm{ev}_t^* (\ch_{\mathrm{odd}}) \in
    \Omega^{\mathrm{even}}(P \grres) 
\end{align*}
on the path spaces based at the identity. They famously fit into the equation
\begin{align*}
    \dd \cs = \ev_1^* \ch - \ev_0^* \ch
\end{align*}
by an application of Stokes' theorem. When we pull back the Chern--Simons forms
to the based loop space in order to get
a form $\cs_{\Omega}$, this identity shows that $\cs_{\Omega}$ is a
transgression of $\ch$ in the path-loop fibration. Using our universal
representatives, we can now associate to a map into $\mathrm U^1$ or
$\grres$, i.e. to a representative for a $K$-theory class, certain differential
forms.  
\begin{dfn}
    \label{dfn:chernmaps}
    Let $M$ be a smooth compact manifold. Then we define the maps
    \begin{align*}
        \Ch\colon \mathrm{Map}(M,\mathrm U^1) &\to \Omega_{\mathrm{cl}}^{\mathrm{odd}}(M) \\ 
        \Ch\colon \mathrm{Map}(M,\grres) &\to
        \Omega_{\mathrm{cl}}^{\mathrm{even}}(M),
    \end{align*}
    given by pullback of the universal Chern forms (Def. \ref{dfn:oddchern} and
    Def. \ref{dfn:evenchernform}). Furthermore, we define the maps 
    \begin{align*}
        \CS\colon \mathrm{Map}(M\times I,\mathrm U^1)\to 
        \Omega^{\mathrm{even}}(M) \\
        \CS\colon \mathrm{Map}(M\times I,\grres)\to 
        \Omega^{\mathrm{odd}}(M) 
    \end{align*}
    given by ``pullback of the universal Chern--Simons forms'' via smooth
    homotopies, i.e. $\CS(H_t) = \int_I H_t^* \ch$.
\end{dfn}
We define a refined notion of homotopy by using these forms,
following \cite[Def. 3.4]{tradler_differential_2015}. It is
designed to retain more information in an equivalence class than just the
isomorphism type of the corresponding bundle. One important feature is that we will have a
well-defined map that assigns to a ``$\CS$-equivalence'' class of maps the pullback of
its universal Chern form, which is only possible up to exact forms for a homotopy
class.
\begin{dfn}\label{dfn:CSequivalence}
    Let $f,g\colon M\to \mathcal{U}$ for $\mathcal{U}\in \left\{ \grres, \mathrm U^1 \right\}$ be smooth
    maps. We say that $f$ and $g$ are Chern--Simons
    homotopic
    (CS-homotopic) if there is a smooth
    homotopy $H_t$ connecting them such that the resulting Chern--Simons form given
    by integrating the universal Chern character
    \begin{align*}
        \CS_{\textup{odd/\textup{even}}}(H) = \int_{I}
        H_t^*(\ch_{\textup{even/\textup{odd}}}) \in
        \Omega^{\textup{odd/\textup{even}}}(M)
    \end{align*}
    is exact.
\end{dfn}
We will also define the blocksum operation, which works in general for operators
on an infinite-dimensional Hilbert space $\h$. It will be used to implement
addition in differential $K$-theory. In order to be explicit, we choose a
specific isomorphism $\rho:\h\to \h\oplus \h$. When a polarization on $\h$ is given,
our isomorphism is designed to respect the grading. 
\begin{dfn}\label{dfn:blocksum}
    Let $\rho\colon \h \to \h\oplus \h$ be the isometric isomorphism
    \begin{align*}
        \rho: e_{2k} \mapsto (e_k,0) \\
        e_{2k+1} \mapsto (0,e_k),
    \end{align*}
    given on a $\N$ or $\Z$-graded orthonormal basis $\left\{ e_i \right\}$. We define
    the corresponding blocksum map
    \begin{align*}
        \boxplus_\rho\colon \mathfrak{gl}(\h)\times \mathfrak{gl}(\h)&\to \mathfrak{gl}(\h) \\
        (A,B)&\mapsto \rho^* (A\oplus B) \rho.
    \end{align*}
\end{dfn}
Note that various subgroups of operators which we consider are preserved by this
construction, most importantly $\ures$ and $\mathrm U^1$. This also induces 
a well-defined operation on $\grres$, where
it corresponds to a direct sum of subspaces: If $W=X(\h_+)$ and $V=Y(\h_+)$
for $X,Y\in \ures$, then
\begin{align*}
    W\boxplus_{\rho} V = (X\boxplus_{\rho} Y)(\h_+) = \rho^*(X\oplus Y)\rho(\h_+) = \rho^*
    (V\oplus W),
\end{align*}
where in the last expression we interpret $V\oplus W$ as a subspace of $\h\oplus
\h$, according to $\rho$.

Via pointwise application, we can now make sense of the blocksum of two maps $f,g$ from
a manifold into the bounded linear operators $\mathfrak{gl}(\h)$. We write
\begin{align*}
    f\boxplus_{\rho} g = \rho^* (f\oplus g) \rho.
\end{align*}
Ultimately, one wants this blocksum operation on maps to not depend on the chosen unitary
isomorphism $\rho$ up to the right equivalence relation. This is easily seen to
be true for homotopy classes of maps by using path-connectedness of the unitary
group $\mathrm U$. The following technical lemma will show the corresponding
statement for the more restricted class of $\CS$-homotopies.
\begin{lem}
    \label{lem:unitaryconjugationisfine}
    Let $f\colon M\to \mathrm U^1$, $g\colon M\to \ures$ and $h\colon M\to
    \grres$ be smooth maps and
    consider $A\in \mathrm U(\h_+)$ and $B\in \mathrm U(\h_+)\times \mathrm
    U(\h_-)\subset \ures$. Then the pairs of maps
    \begin{enumerate}[(i)]
        \item $A f A^*\colon M\to \mathrm U^1 \qquad\; \text{and} \qquad f\colon
            M\to \mathrm U^1$ 
        \item $Bg B^* \colon M\to\ures\quad\;\,\; \text{and} \qquad g\colon M\to
            \ures $
        \item $B h \colon M\to\grres\qquad \text{and} \qquad h\colon M\to \grres$
    \end{enumerate}
    are $\CS$-homotopic\footnote{On $\ures$, $\CS$-equivalence is defined with
        respect to the universal Chern character that one gets from pulling back
        $\ch_{\mathrm{even}}$ via the projection $\ures\to\grres$.}, i.e.\
        conjugation by a fixed such matrix does not change
    the Chern--Simons equivalence class. In particular, for any other unitary
    isomorphism $\rho'\colon \h \to \h\oplus \h$ (respecting the grading in the
    polarized case), we have that $f\boxplus_{\rho}
    g$ and $f\boxplus_{\rho'} g$ are $\CS$-homotopic.
    \begin{proof}
        For the first case, choose a smooth path $A_t$ from $A_0 = 1$ to $A_1=
        A$. Then there is a smooth universal homotopy 
        \begin{align*}
            H_t\colon \mathrm U^1 \times I &\to \mathrm U^1 \\ (X,t)&\mapsto A_t
            X A_t^*,
        \end{align*}
        which yields a homotopy as stated for any $f\colon X\to \mathrm U^1$ by
        composition. We need to show that its $\CS$-form is exact. We have
        \begin{align*}
            \dd \CS_{2k}(H_t) = \dd \int_{I}H_t^*\ch_{2k+1} = H_1^* \ch_{2k+1} -H_0^*
            \ch_{2k+1} = \\-\left(\frac{1}{2\pi i}\right)^{k+1}\frac{(k)!}{(2k+1)!} (\tr
            (A X^* \dd X A^*)^{2k+1} - \tr (X^* \dd X)^{2k+1})=0. 
        \end{align*}
        Since the positive even cohomology of $\mathrm U^1$ vanishes, this
        implies that the Chern--Simons forms for $k>0$ are exact. For $k=0$, we
        make a direct calculation. Note that the differential
        of $H_t$ splits according to the splitting of the tangent space of
        $\mathrm U^1\times I$ into a 
        sum of a space part with a time derivative. Our notation for the
        space derivative is $\dd
        H_t$, while we denote the time derivative by $\dot H_t$. We have
        \begin{align*}
            \dot H_t &= \dot A_t X A_t^* - A_t X A_t^* \dot
            A_t A_t^* \\ 
            \dd H_t &= A_t \dd X A_t^*.
        \end{align*}
        We need to calculate
        \begin{align*}
            \CS_{0}(H_t) = \int_{I}H_t^*\ch_{1} = \left(\frac{i}{2\pi}\right)\int_I
            \iota_{\partial_t} ( H_t^* (\tr (\omega_{\mathrm{MC}}))).
        \end{align*}
        The integrand yields
        \begin{align*}
            \iota_{\partial_t} ( H_t^* (\tr (\omega_{\mathrm{MC}}))) &= \tr (A_t X^* A_t^* (\dot A_t X A_t^* -
            A_t X A_t^* \dot A_t A_t^*)) \\&= \tr (X^* A_t^* \dot A_t X - \dot A_t
            A_t^*) = 0
        \end{align*}
        for all $t$ and therefore, $\CS_0$ vanishes.

        For the second case, we choose again a smooth path from $1$ to
        $B$ in order to define a homotopy
        $H_t$ starting at $H_0=C_B$ and ending at $H_1=\id_{\ures}$, where $C_B$ denotes
        conjugation by $B$. By the vanishing
        of $H^{\mathrm{odd}}(\grres)$, it is enough to show that the $\CS$-form
        is closed, i.e. that $H_0$ and $H_1$ have the same Chern form. We argue
        as follows: The
        projection $\ures\to \grres$ sends a matrix $X$ to the projection
        $X\pi_+ X^*$. Conjugating $X$ by $B$ yields
        \begin{align*}
            BXB^* \pi_+ B X^* B^* = B X \pi_+ X^* B^*.
        \end{align*}
        and therefore using the invariance of $\ch$, the conjugation map $C_B\colon
        \ures\to \ures$ pulls back the universal Chern form to itself, i.e.
        $C_B^* \ch = \ch = \id_{\ures}^*\ch$. The third case follows by the same
        argument, using the invariance of $\ch$ one more time. 
        
        The independence of the blocksum up to $\CS$-equivalence is now easily
        deduced, since 
        \begin{align*}
            {\rho'}^* \rho(f\boxplus_{\rho}g)\rho^* {\rho'} = {\rho'}^*
            \rho\rho^* (f\oplus g) \rho \rho^* {\rho'} = f\boxplus_{\rho'} g
        \end{align*}
        and therefore the two blocksums defined by $\rho$ and $\rho'$ just differ by a conjugation with the
        unitary matrix $\rho'\rho^*$ on $\h$, which in the polarized case
        repects the grading.
    \end{proof}
\end{lem}
\begin{rem}
    By the preceding Lemma, it is now safe to suppress $\rho$ in our
    notation. For two elements in
    the restricted unitary group, which are by definition 2 by 2 block
    operators, we write
    \begin{align}
        \label{eq:blocksumfg}
        f\boxplus g = \begin{pmatrix}
            f_{++} & 0 & f_{-+} & 0 \\ 0 & g_{++} & 0 & g_{-+} \\ f_{+-} & 0 &
            f_{--} & 0 \\ 0 & g_{+-} & 0 & g_{--}
        \end{pmatrix}.
    \end{align}
\end{rem}
\begin{prop}
    \label{prop:commutativity}
    Let $f,g,h\colon M\to \mathcal{U}$ be smooth maps for $\mathcal{U}
    \in \left\{ \grres, \mathrm U^1 \right\}$. Then, the operation induced by blocksum is
    commutative and associative up to $CS$-homotopy, i.e. we have
    \begin{align*}
        f \boxplus g \sim_{\CS} g \boxplus f \quad \mathrm{and} \quad f\boxplus (g
        \boxplus h) \sim_{\CS}(f\boxplus g )\boxplus h.
    \end{align*}
    \begin{proof}
        This is just a consequence of Lemma \ref{lem:unitaryconjugationisfine},
        since the difference in each case is just a permutation of the basis.
        For commutativity in the case of $\mathrm U^1$, one sees that
        for $U=
        \begin{pmatrix}
            0&1\\1&0
        \end{pmatrix}\in \mathrm U(\h\oplus \h)$, one has 
        \begin{align*}
            g\boxplus f = \rho^*U\rho (f\boxplus g) \rho^* U \rho.
        \end{align*}
        For the even case, acting by the same matrix $\rho^* U \rho$ on $\h_+$ and $\h_-$
        seperately does the job. For associativity, one has that
        \begin{align*}
            (f\boxplus g)\boxplus h &= \rho^* (\rho^*\times \id) \begin{pmatrix}
                f & &\\ &g&\\&&h
            \end{pmatrix} (\rho\times \id) \rho \\ & = \rho^* (\rho^*\times \id)
            (\id\times \rho)\rho (f\boxplus (g\boxplus h))\rho^* (\id\times\rho^*)(\rho\times \id) \rho
        \end{align*}
        in the $\mathrm U^1$ case and acting by the same matrix on $\h_+$ and $\h_-$
        separately does the job in the $\grres$ case.
    \end{proof}
\end{prop}
We will now discuss the involution on $\ures$ that will implement inversion in
differential $K$-theory. Let $U$ be the unitary transformation that flips the
role of $\h_+$ and $\h_-$, given by the matrix $U = \begin{pmatrix} 0 & 1 \\ 1 & 0
\end{pmatrix}$. In a $\Z$-basis $\{e_i\}$ adapted to the polarization, the map
$U$ sends $e_i$ to $e_{-i-1}$.
\begin{dfn}
    \label{dfn:flipmap}
    We define the polarization flip map $\flip\colon \ures\to\ures$ to be
    conjugation with $U$. On the space of smooth maps from a manifold to
    $\ures$, this induces the operation 
    \begin{align*}
        f\mapsto \mathrm{flip} (f) = \flip \circ f = UfU.
    \end{align*}
    Explicitly, we have
    \begin{align*}
        \mathrm{flip}(f)(x) = \begin{pmatrix}
            f_{--}(x) & f_{+-}(x) \\ f_{-+}(x) & f_{++}(x)
        \end{pmatrix}, \;\textrm{when} \quad f(x) = \begin{pmatrix}
            f_{++}(x) & f_{-+}(x) \\ f_{+-}(x) & f_{--}(x)
        \end{pmatrix}.
    \end{align*}
\end{dfn}
Note that there is an induced flip map on the restricted Grassmannian, which
corresponds to taking the orthogonal complement of a subspace and then changing
the polarization. We have
\begin{align*}
    W= X(\h_+) \mapsto \flip(X)(\h_+) = U X U(\h_+) = U X(\h_-) = U(W^{\perp}) =
    \flip(W),
\end{align*}
which also extends to maps $M\to\grres$ via composition. One furthermore sees that pullback
by flip preserves left invariance of forms: If
$L_Y^*\eta = \eta$, we have that
\begin{align*}
    L_Y^*(\mathrm{flip}^* \eta) =(\mathrm{flip}\circ L_Y)^* (L_{\mathrm{flip}(Y^{-1})}^* \eta)
    = \mathrm{flip}^* \eta,
\end{align*}
since $\mathrm{flip}$ is a group homomorphism on $\ures$. The following
proposition shows compatibility of the inversion and addition operations on the
classifying spaces with the Chern and Chern--Simons forms.
\begin{prop}
    \label{prop:CSforms}
    Consider smooth maps $f,g\colon M\to \mathcal{U}$ and $\mathcal{U}\in
    \left\{ \grres, \mathrm U^1 \right\}$. Then:
    \begin{enumerate}[(i)]
        \item The maps $\Ch$ and $\CS$ are monoid morphism, i.e. $\Ch(f\boxplus
            g) = \Ch(f) + \Ch(g)$ and $\CS(H_t \boxplus G_t) = \CS(H_t) + \CS(G_t)$.
        \item $\CS(H_t * G_t) = \CS(H_t) + \CS(G_t)$ (composition of
            homotopies).
        \item $\Ch_{\mathrm{even}}(\flip (f)) = -\Ch_{\mathrm{even}}(f)$, and
            $\Ch_{\mathrm{odd}}(f^*) = -\Ch_{\mathrm{odd}}(f)$.
        \item $\CS_{\textup{odd}}(\flip (H_t)) = - \CS_{\textup{odd}}(H_t)$,
            and $\CS_{\textup{even}}(H^*_t) = -\CS_{\textup{even}}(H_t)$.
    \end{enumerate}
    \begin{proof}
        The monoid morphism property follows directly from the additivity of the trace
        under block sum and linearity of the integral, and the additivity under
        composition follows from additivity of the integral under partition of
        the intervall.
       
        We check the third identity directly on $\grres$ and $\mathrm U^1$. For
        $\grres$, we need to compute the pullback of the curvature $\flip^* \Omega$, and it
        suffices to do this in the tangent space at $\h_+$ by left
        invariance. Take $X=\begin{pmatrix}
            0&X_{-+}\\X_{+-}&0
        \end{pmatrix}$ and $Y=\begin{pmatrix}
            0&Y_{-+}\\Y_{+-}&0
        \end{pmatrix}\in T_{\h_+}\grres$. We have
        \begin{align*}
            (\flip^*\Omega)_{\h_+}(X,Y) = \Omega\left(\begin{pmatrix}
                0& X_{+-} \\ X_{-+} &0 
            \end{pmatrix}, \begin{pmatrix} 
                0&Y_{+-} \\Y_{-+} & 0 \end{pmatrix}\right) =
            Y_{+-}X_{-+} - X_{+-}Y_{-+}
        \end{align*} 
        Therefore, we see that $(\flip^* (\tr\,\Omega^k))_{\h_+} (X^1 ,\dots, X^{2k})$
        equals
        \begin{align*}
            \frac{1}{2^{k}}\sum_{\sigma\in S_{2k}} \mathrm{sgn}(\sigma)\tr
            \,\Omega(\flip X^{\sigma(1)},\flip X^{\sigma(2)})
            \cdots \Omega(\flip X^{\sigma(2k-1)}, \flip X^{\sigma(2k)})\\ =  
            \frac{1}{2^{k}}\sum_{\sigma\in S_{2k}}\mathrm{sgn}(\sigma)\tr \,(X^{\sigma(2)}_{+-}
            X^{\sigma(1)}_{-+} -
            X^{\sigma(1)}_{+-} X^{\sigma(2)}_{-+}) \cdots (X^{\sigma(2k)}_{+-}
            X^{\sigma(2k-1)}_{-+} - X^{\sigma(2k-1)}_{+-} X^{\sigma(2k)}_{-+}).
        \end{align*}
        This is just a big sum of products of $2k$-operators with many redundant
        terms. One sees that it is equal to
        \begin{align*}
            &\sum_{\sigma\in S_{2k}}\mathrm{sgn}(\sigma)\tr\,X^{\sigma(2)}_{+-}X^{\sigma(1)}_{-+}
                X^{\sigma(4)}_{+-}X^{\sigma(3)}_{-+}\cdots
                X^{\sigma(2k)}_{+-}X^{\sigma(2k-1)}_{-+} \\ &= 
            \sum_{\sigma\in S_{2k}}\mathrm{sgn}(\sigma)\tr\,X^{\sigma(1)}_{-+}
                X^{\sigma(4)}_{+-}X^{\sigma(3)}_{-+}\cdots
                X^{\sigma(2k)}_{+-}X^{\sigma(2k-1)}_{-+}X^{\sigma(2)}_{+-} \\&=
            -\sum_{\sigma\in S_{2k}}\mathrm{sgn}(\sigma)\tr\,X^{\sigma(2)}_{-+}
                X^{\sigma(1)}_{+-}X^{\sigma(4)}_{-+}X^{\sigma(3)}_{+-}\cdots
                X^{\sigma(2k)}_{-+}X^{\sigma(2k-1)}_{+-} \\&=
                -(\tr\, \Omega^k)_{\h_+}(X^1,\dots,X^{2k})
        \end{align*}
        The first equality is cyclic invariance of the trace, the second one
        comes applying a cyclic permutation with $2k$-elements. Finally, the
        last equality comes from going through the same calculation without
        applying the flip map. This proves (iii) for the odd case.

        For even case, notice that from $0 = \dd (f f^*) = \dd f f^* + f \dd
        f^*$ it follows that $\tr (f \dd f^*)^{2k-1} = - \tr (\dd f f^*)^{2k-1} =- \tr (
        f^*\dd f)^{2k-1}$. Therefore pulling back $\ch_{\mathrm{odd}}$ via the
        adjoint-operator map $\mathrm U^1\to \mathrm U^1$ gives the desired
        minus sign. Part (iv) easily follows from part (iii) by definition of
        the Chern--Simons form.
    \end{proof}
\end{prop}

\section{Geometric structure maps in the $K$-theory spectrum}
\label{sec:smoothper}

The goal of this section is to recall explicit homotopy equivalences $\Omega
\grres \to \mathrm{U}^1$ and $\Omega \mathrm{U}^1 \to \grres$ which are
compatible with our Chern and Chern--Simons forms. We will consider the even case
first.

In Section \ref{sec:restricted}, we constructed a smooth model for the universal
$\mathrm{GL}^1$ principal bundle as $\stres\to \grres^0$, together with an
explicit connection form $\Theta$. Consider the loop space
$\Omega\grres$, based at $\h_+\in \grres^0$. Parallel transport via the
connection $\Theta$ gives rise to the holonomy map, which assigns to such a loop
the fiber coordinate of the endpoint of the horizontal lift of this loop,
starting at $w_0 = \begin{pmatrix}
    1\\0
\end{pmatrix}\in \stres$. After applying the homotopy equivalence
\begin{align*}
    \mathrm{GL}^1 &\to \mathrm{U}^1 \\ T &\mapsto T |T|^{-1},
\end{align*}
we have constructed a map $h_{\mathrm{even}}\colon
\Omega\grres\to \mathrm U^1$. Since it implements
holonomy in the fibration $\mathrm{U} \to E \mathrm{U} \to B\mathrm{U}$, it is
clear that this is a homotopy equivalence. It remains to check compatibility
with the Chern and Chern--Simons forms.

It is well known that the transgression of $[\ch_{\mathrm{even}}]\in
H^{\mathrm{even}}(B\mathrm U; \R)$ in the
universal fibration $\mathrm U\to E \mathrm U\to B\mathrm U$ is the class
$[\ch_{\mathrm{odd}}]\in H^{\mathrm{odd}}(U;\R)$. We can actually recover this
fact from the universal bundle which we constructed in Section
\ref{sec:restricted} by a direct calculation.
\begin{lem}
    \label{lem:transgression}
    The transgression map $T\colon H^k(\grres^0;\R)\to H^{k-1}(\stres\colon \R)$ in the universal
    $\mathrm{GL}^1$-fibration $\stres\to \grres^0$ maps the even Chern character
    to the odd one, i.e. $T([\ch_{2k}]) = [\ch_{2k-1}]$.
    \begin{proof}
        We use the connection $\Theta$ as constructed in Prop.
        \ref{prop:univconnection}. Then the Chern character is given by an
        invariant polynomial, evaluated at the curvature $\Omega$. For this situation,
        Chern and Simons \cite[Sec. 3]{chern1974characteristic} gave a formula for the
        transgression. Define $\varphi_t = t \Omega + 1/2 (t^2
        -t) [\Theta,\Theta]$. Then one easily checks that the form
        \begin{align*}
            \eta = \left(\frac{i}{2\pi}\right)^k \frac{1}{k!}\int_{0}^1 k \tr
            (\Theta \wedge \varphi_t^{k-1}) \dd t
        \end{align*}
        satisfies the two identities $\dd \eta = \pi^*\ch_{2k}$ and $i^*\eta =
        \ch_{2k-1}$. Therefore $\ch_{2k-1}$ represents a transgression of
        $\ch_{2k}$.
    \end{proof}
\end{lem}
From this result we can deduce the needed compatibility between the Chern and
Chern--Simons forms.
\begin{prop}
    \label{prop:compatibilityeven}
    Let $H\colon M\times I \to \grres^0$ be a smooth homotopy, starting and
    ending with $H_0=H_1=\mathrm{const}_{\h_+}$ with adjoint map $\hat
    H\colon M\to \Omega\grres^0$. Let $h_{\mathrm{even}}\colon
    \Omega\grres^0\to \mathrm{GL}^1\to\mathrm{U}^1$ be the holonomy map composed with the
    homotopy equivalence $\mathrm{GL}^1\to \mathrm{U}^1$ given by $X\mapsto X
    |X|^{-1}$. Then, on the level of differential forms, we have the congruence
    modulo exact forms
    \begin{align*}
        \int_I H_t^*\ch_{\mathrm{even}} = (h_{\mathrm{even}}\circ \hat H)^*\ch_{\mathrm{odd}}
        + \mathrm{exact}.
    \end{align*}
    \begin{proof}
        We have a diagram of fibrations
        \[\begin{tikzcd}
                \mathrm{GL}^1 \arrow[r]& \stres \arrow[r]& \grres^0\\
                \Omega \grres^0 \arrow[r]\arrow[u,"\mathrm{hol_\Omega}"] & P \grres^0 \arrow[r]\arrow[u,
                "\mathrm{hol}"] & \grres^0\arrow[u, equal],
        \end{tikzcd}\]
        where the vertical maps are induced by holonomy. Now since
        transgression commutes with maps of fibrations, we see that 
        $\int_{I} \mathrm{ev}_t^* [\ch_{\mathrm{even}}]=
        \mathrm{hol}_{\Omega}^*[\ch_{\mathrm{odd}}] = h_{\mathrm{even}}^*[\ch_{\mathrm{odd}}]$
        as cohomology classes on $\Omega \grres^0$.

        We now pull back the equation via $\hat H$ and get 
        \begin{align*}
            (h_{\mathrm{even}}\circ \hat H)^*[\ch_{\mathrm{odd}}] = \hat H^* \int_{I}
            \mathrm{ev}_t^* [\ch_{\mathrm{even}}]= \int_{I} (\hat H\times
            \id_I)^* \mathrm{ev}_t^* [\ch_{\mathrm{even}}] = \int_{I}
            H_t^*[\ch_{\mathrm{even}}] 
        \end{align*}
        Since the domain now is a finite-dimensional manifold, we can smoothly
        approximate the maps up to homotopy and see that the claimed equality is
        true on the level of differential forms, up to exact forms.
    \end{proof}
\end{prop}
\begin{rem}
    \label{rem:avoidsmoothness}
    This argument (and also the one in the odd case given below) avoids the
    discussion of $h_{\mathrm{even}}\colon \Omega\grres^0\to \mathrm
    U^1$ as a smooth map between infinite-dimensional manifolds. It would be interesting to
    compute directly the derivative of $h_{\mathrm{even}}$ and try to pull back
    $\ch_{\mathrm{odd}}$ as a differential form.
\end{rem}
The construction of the map for the odd case has appeared in
\cite[Appendix 2]{Carey2000}, based on the ideas of \cite[Sec. 6.3]{pressley_loop_1988}.
We choose a concrete model for the generic polarized Hilbert space
that has been used before. Let $\h^{\infty} = L^2(S^1, \h) = L^2(S^1)
\hat\otimes\h$ be the space of 
$L^2$-functions on the circle to an infinite-dimensional separable complex
Hilbert space $\h$ with fixed basis $\left\{ e_i \right\}_{i\geq 0}$. There is a
natural $\Z_2$-grading given by
the positive resp.\ negative exponent part of the Fourier decomposition, i.e.
\begin{align*}
    \h^{\infty}_+ = \left\{ f\in \h^{\infty} \mid f = \sum_{k\geq 0} f_k
    z^k, f_k\in \h\right\}\\ 
    \h^{\infty}_- = \left\{ f\in \h^{\infty} \mid f = \sum_{k< 0} f_k
    z^k, f_k\in \h\right\},
\end{align*}
where $z = \exp(i \theta)$. We consider the multiplication operator map
\begin{align*}
    h_{\mathrm{odd}}\colon\Omega \mathrm U^1(\h)&\to \mathrm U(\h^{\infty})\\ \gamma & \mapsto M_\gamma,
\end{align*}
where $(M_\gamma f)(\theta) = \gamma(\theta) f(\theta)$. We have the following 
\begin{lem}
    \label{lem:multiplication}
    The map $h_{\mathrm{odd}}$ has image in the restricted unitary group. Furthermore, as a map
    to $\ures$, $h_{\mathrm{odd}}$ is a homotopy equivalence.
    \begin{proof}
        We will rely on the corresponding statements for the finite-dimensional
        version of said map which were proved by
        Pressley and Segal in \cite[Ch. 8]{pressley_loop_1988}.
        
        First, we have that the analogously defined map $h_{\mathrm{odd}}^n\colon \Omega
        \mathrm U(n)\to \mathrm U(\h^{(n)})$ has image contained in
        $\ures(\h^{(n)})\subset \ures(\h^{\infty})$, where $\h^{(n)}$ is the
        finite-dimensional version of our space $\h^{\infty}$, i.e. $\h^{(n)} =
        L^2(S^1,\C^n)$, and the inclusion $\C^n\hookrightarrow \h$ is via the
        first $n$ basis vectors
        $\left\{ e_i \right\}_{i\in \N}$. This is a consequence of the decay condition on the Fourier
        coefficients of the loop $\gamma$, using the boundedness of its first derivative. 
        Furthermore, we have that this map is $(2n-2)$-connected. 

        In order to conclude the corresponding statements for our stabilized
        version of the map, we note that the restriction of $h_{\mathrm{odd}}$ to $\Omega
        \mathrm U(n)$
        has image contained in $\ures$ by the previous paragraph. The union of
        these loop spaces is the loop space $\Omega \mathrm U$ of the stable unitary
        group, and this still maps to $\ures$. But $\Omega \mathrm U$ is a dense
        subspace of $\Omega \mathrm U^1$ and therefore the first claim follows, since
        $\ures$ is a complete Riemannian manifold.

        Since the connectivity of these maps increases, we can use a similar
        argument for $h_{\mathrm{odd}}$ being a homotopy equivalence. It is enough to show that
        it induces an isomorphism on all homotopy groups. We have the
        commutative diagram 
        \[\begin{tikzcd}
                \Omega \mathrm U(n) \arrow[r, "h_{\mathrm{odd}}^{n}"]\arrow[d, "i"] &
                \ures(\h^{(n)})\arrow[d, "j"]\\
                \Omega \mathrm U^1(\h) \arrow[r, "h_{\mathrm{odd}}"] & \ures(\h^{\infty}),
        \end{tikzcd}\]
        where both horizontal maps come from the inclusion of $\C^n$ into
        $\h$ and then filling up with the identity matrix. The maps $h_{\mathrm{odd}}^{n}$
        and $i$ are $2n-2$-connected, while the map $j$ is a homotopy
        equivalence and therefore $h_\mathrm{odd}$ is also $2n-2$-connected, for any
        $n$.
    \end{proof}
\end{lem}
The map $h_{\mathrm{odd}}$ realizes the inverse of the Bott periodicity map as a homomorphism
of infinite-dimensional Lie groups. We can append the projection $\ures\to\grres$
in order to get our desired periodicity map for our Grassmannian model, and we
will also denote this map by $h_{\mathrm{odd}}$. 
\begin{prop}
    \label{prop:compatibilityodd}
    Let $H\colon M\times I \to \mathrm U^1$ be a smooth homotopy, starting and
    ending with $H_0=H_1=\mathrm{const}_{\mathrm{\id}}$ with adjoint map $\hat
    H\colon M\to \Omega\mathrm U^1$. Let $h_{\mathrm{odd}}\colon
    \Omega \mathrm U^1 \to \ures \to \grres$ be the assignment of the
    corresponding multiplication operator, composed with the homotopy
    equivalence given by the projection. Then on the level of differential
    forms, we have the congruence modulo exact forms
    \begin{align*}
        \int_I H_t^*\ch_{\mathrm{odd}} = (h_{\mathrm{odd}}\circ \hat H)^*\ch_{\mathrm{even}}
        + \mathrm{exact}.
    \end{align*}
    \begin{proof}
        In the proof of Prop. \ref{prop:compatibilityeven}, we have seen that
        the even Chern character transgresses to the odd Chern character as
        cohomology classes in the path loop fibration over $\grres^0$. Now since
        the Chern character is compatible with Bott periodicity, if we apply
        transgression again in the path loop fibration over $\mathrm U^1$, we
        must get $T^2([\ch_{2k}]))= [\ch_{2k-2}]$ after identifying $\Omega^2
        \grres^0$ with $\grres$ explicitely via $h_{\mathrm{odd}}\circ (\Omega
        h_{\mathrm{even}})$. This allows
        us to make the calculation
        \begin{align*}
            (h_{\mathrm{odd}}\circ \Omega h_{\mathrm{even}})^*[\ch_{\mathrm{even}}] = T (T
            ([\ch_{\mathrm{even}]})) = T
            (h_{\mathrm{even}}^*[\ch_{\mathrm{odd}}]) = (\Omega
            h_{\mathrm{even}})^* T ([\ch_{\mathrm{odd}}]),
        \end{align*}
        where the first equality is Bott periodicity, the second one is from
        (the proof of) Prop. \ref{prop:compatibilityeven}, and the last one
        follows from the naturality of transgression in the diagram of
        fibrations
        \[\begin{tikzcd}
                \Omega\mathrm{U}^1 \arrow[r]& P \mathrm U^1 \arrow[r]& \mathrm
                U^1\\
                \Omega^2 \grres^0 \arrow[r]\arrow[u,"\mathrm{\Omega h_{\mathrm{even}}}"] & P
                \Omega\grres^0 \arrow[r]\arrow[u,
                "P h_{\mathrm{even}}"] & \Omega\grres^0\arrow[u, "h_{\mathrm{even}}"].
        \end{tikzcd}\]
        Since $\Omega h_{\mathrm{even}}$ is a homotopy equivalence, we get that $T
        [\ch_{\mathrm{odd}}] = h_{\mathrm{odd}}^* [\ch_{\mathrm{even}}]$. Pulling back by
        $H$ and using the same argument as in the even case now gives the claim.
    \end{proof}
\end{prop}
\begin{prop}
    \label{prop:spectrum}
    The sequence of pointed spaces and pointed maps $(E_n,h_n)$, $n\in \Z$ given by 
    \begin{align*}
        E_{2n} &= \grres \quad\mathrm{and}\quad E_{2n+1} = \mathrm U^1 \\
        h_{2n} &= h_{\mathrm{even}} \quad\mathrm{and}\quad h_{2n+1} =
        h_{\mathrm{odd}} 
    \end{align*}
    defines an $\Omega$-spectrum that represents complex topological $K$-theory.
    Furthermore, addition in
    $K$-theory is implemented by the blocksum operation on both $\grres$ and
    $\mathrm U^1$.
    \begin{proof}
        Since the structure maps are homotopy equivalences, the first part of
        the theorem follows. For the second part, we have to prove that
        blocksum is homotopic to composition of loops, i.e. that the squares
        \[\begin{tikzcd}
            \Omega \mathrm U^1 \times \Omega \mathrm U^1\arrow[d,
            "h_{\mathrm{odd}}\times h_{\mathrm{odd}}"]\arrow[r, "*"] &
            \Omega \mathrm U^1\arrow[d, "h_{\mathrm{odd}}"] &
            &\Omega \grres \times \Omega \grres \arrow[d,
            "h_{\mathrm{even}}\times h_{\mathrm{even}}"]\arrow[r, "*"] &
            \Omega \grres\arrow[d, "h_{\mathrm{even}}"]
            \\ 
            \grres\times\grres \arrow[r, "\boxplus"]& \grres&
            &\mathrm U^1 \times \mathrm U^1 \arrow[r, "\boxplus"]& \mathrm U^1.
        \end{tikzcd}\]
        commute up to homotopy, where the star denotes loop composition. 
        
        For the left one, recall that
        $h_{\mathrm{odd}}$ assigns to a loop the corresponding multiplication
        operator in $\ures$ and then projects to $\grres$. Since the projection
        $\ures\to\grres$ is a map of $H$-spaces for the blocksum, we can work on
        $\ures$. It is clear that $h_{\mathrm{odd}}$ respects the alternative
        $H$-space structures on target and domain:
        the pointwise multiplication of two loops maps to the product of their
        operators. Now by the usual Eckmann--Hilton argument, pointwise
        multiplication of loops is homotopic to loop concatenation. On the
        other hand, if $C_t$ denotes the grading preserving rotation 
        \begin{align*}
            C_t = \begin{pmatrix}
                \cos(t)&\sin(t) & 0 & 0 \\ -\sin(t) & \cos(t)&0&0 \\
                0&0&\cos(t)&\sin(t)\\ 0&0&-\sin(t)&\cos(t)
            \end{pmatrix},
        \end{align*}
        we have the homotopy
        \begin{align*}
            \ures\times\ures\times I &\to \ures \\ (A,B,t)&\mapsto (A\boxplus 1)
            C_t (1\boxplus B) C_t^*,
        \end{align*}
        which shows that $A\boxplus B \sim AB \boxplus 1$. We will finish the proof by
        showing that the map $f\colon A\mapsto A\boxplus 1$ is homotopic to the
        identity. It is enough to show this over each path component of $\ures$,
        so let $\ures^0$ denote the identity component. We have homotopy
        equivalences $B\mathrm U \stackrel{\varphi}{\to}  \ures^0 \stackrel{\psi}{\to}
        \mathrm{Fred}^0$ and the second map
        is given explicitly by the projection to the $++$-component (recall 
        (\ref{eq:htpyeq})). Notice that this respects taking the blocksums of
        operators, i.e. $f\circ \psi = \psi\circ f$.
        
        We have $B\mathrm U = \colim B\mathrm U(n)$. Since each of the
        inclusions in this colimit induces a
        surjection in $K$-theory, there is no $\lim^1$-term, and we have
        \[\begin{tikzcd}[contains/.style = {draw=none,"\in" description,sloped}]
                \left[\ures^0,\ures^0 \right] \arrow[r, "\cong"]& \left[B\mathrm U,
                \mathrm{Fred}^0\right] \arrow[r, "\cong"]& \left[\colim B\mathrm U(n),
                \mathrm{Fred}^0\right]\arrow[r, "\cong"]& \lim [B\mathrm U(n),
                \mathrm{Fred}^0 ] \\ \left[\mathrm{f}\right]\arrow[u,
                contains]\arrow[r, mapsto] & \left[\psi\circ f\circ \varphi\right]\arrow[u,
                contains]\arrow[r, equals] &\left[f\circ \psi\circ \varphi\right]\arrow[u,
                contains]\arrow[r, mapsto] &(\left[\restr{f\circ\psi\circ
                \varphi}{B\mathrm U(n)}\right])\arrow[u,
                contains]
        \end{tikzcd}\]
        But $B\mathrm U(n)$ is compact, and by the Atiyah--Jänich index map,
        $[B\mathrm U(n),\mathrm{Fred}^0]\cong \widetilde K(B\mathrm U(n))$.
        Since blocksum with the identity does not change the kernel or cokernel
        of an operator, we have that $(\left[\restr{f\circ\psi\circ
        \varphi}{B\mathrm U(n)}\right]) = (\left[\restr{\psi\circ
        \varphi}{B\mathrm U(n)}\right])$, and tracing back through the
        isomorphisms, we see that $[f] = [\mathrm{id}]$.
        
        The argument for the second square in the above diagram is similar and
        follows from the fact that the holonomy map $h_{\mathrm{even}}$ takes
        compositions of loops to products of operators, and that $\mathrm
        U^1\simeq \mathrm U$ is filtered by the compact manifolds $\mathrm U(n)$
    \end{proof}
\end{prop}

\section{The Differential $K$-theory groups}
\label{sec:differential}

Let us recall the axiomatic definition of
differential extensions according to \cite[Def. 1.1]{bunke2010uniqueness}, adapted to the
special case of $K$-theory. Let $V = K^*(*)\otimes_Z\R = \R[u,u^{-1}]$ be the
coefficients of complex $K$-theory with the Bott element of degree $2$, and
let $\Omega^*(M;V) =
\mathscr{C}^{\infty}(M,\Lambda^* T^* M \otimes_\R V)$ denote the $V$-valued
differential forms on $M$, where the degree is induced by the sum of the degrees
as a differential form and as an element of $V$. There
is also a version of cohomology with coefficients in the graded vector space
$V$, which we will denote by $H^*(M;V)$. In the following, we will consider all
$\Z$-graded vector spaces that arise as $\Z_2$-graded by restricting to even or
odd degrees. Denote by $\ch\colon K^*(M)\to H^*(M;V)$ the
topological Chern character.
\begin{dfn}
    \label{dfn:differentialextensions}
    A differential extension of $K$-theory is a contravariant functor from the
    category of smooth manifolds to $\Z_2$-graded abelian groups, together with
    natural transformations
    \begin{enumerate}[(i)]
        \item $R\colon \hat K^*(M) \to \Omega^{*}_{d=0}(M;V)$, called the curvature,
        \item $I\colon \hat K^*(M)\to K^*(M)$, called the underlying class,
        \item $a\colon \Omega^{*-1}(M;V)/ \mathrm{im}(d)\to \hat K^*(M)$, called
            the action of forms,
    \end{enumerate}
    such that
    \begin{enumerate}[(i)]
        \item the diagram
        \[\begin{tikzcd}
            \hat K^*(M) \arrow[d, "R"]\arrow[r, "I"] & K^*(M)\arrow[d, "\ch"] \\ 
            \Omega^*_{d=0}(M;V) \arrow[r, "\mathrm{Rham}"]& H^*(M;V).
        \end{tikzcd}\]
        commutes,
        \item $R\circ a=d$, so the action map is a lift of the exterior
            derivative, and
        \item we have the exact sequence
        \[\begin{tikzcd}
                K^{*-1}(M) \arrow[r, "\mathrm{Ch}"] & \Omega^{*-1}(M;V) / \mathrm{im}(d)
                \arrow[r, "a"] & \hat K^*(M) \arrow[r, "I"]& K^*(M) \arrow[r] &0.
        \end{tikzcd}\]
    \end{enumerate}
\end{dfn}
We will now give a concrete implementation of such a differential refinement via
smooth classifying spaces. The goal of this section is to define the underlying
group-valued functors $\hat K^*$.
\begin{dfn}\label{def:Ktheory}
    Let $M$ be a smooth manifold and $\h\cong \h_+\oplus\h_-$ be an
    infinite-dimensional $\Z_2$-graded separable complex Hilbert space with both
    $\h_+$ and $\h_-$ infinite-dimensional. Define
    the underlying sets of the odd and even differential $K$-theory groups via
    \begin{align*}
        \hat K^0(M) = \textrm{Map}(M, \grres(\h)) /\sim\\
        \hat K^{1}(M) = \textrm{Map}(M, \mathrm U^1(\h_+)) /\sim. 
    \end{align*}
    The equivalence relation is induced by Chern--Simons
    homotopy equivalence (see Def. \ref{dfn:CSequivalence}), together with a stabilization step that identifies any map
    $f$ with $f\boxplus \mathrm{const}_{*}$, where $*$ is the basepoint, i.e. the
    subspace $\h_+$ or the identity.
\end{dfn}
\begin{lem} \label{lem:boxplusgroup}
    The operation $\boxplus$ induces an abelian group structure on $\hat K$. The
    neutral elements are given by the equivalence class of the constant map to
    the basepoint, and inversion is given by $f\mapsto f^*$ and $f\mapsto
    \flip(f)$ in the odd/even case respectively.
    \begin{proof}
        We need to check well-definedness. If $f_0\sim_{\CS} f_1$ and $g_0
        \sim_{\CS} g_1$, then we need to show that
        $f_0 \boxplus g_0 \sim_{\CS} f_1 \boxplus g_1$. This is achieved by the
        homotopy $f_t \boxplus g_t$, which is again a CS-homotopy by Prop.
        \ref{prop:CSforms}. Furthermore, the matrices
        \begin{align*}
                (f\boxplus 1) \boxplus (g\boxplus 1) \qquad \mathrm{and} \qquad
                (f\boxplus g)\boxplus 1 = (f\boxplus g)\boxplus (1 \boxplus 1)
        \end{align*}
        are CS-equivalent by Prop. \ref{prop:commutativity}, and so stabilization is
        also fine. Commutativity and associativity are also proven in Prop.
        \ref{prop:commutativity}. That blocksumming with $\mathrm{const}_{*}$ is
        the identity is built into the definition of our equivalence relation.
        It remains to show that inversion is given by the proposed operations.

        Start with the even case and consider the rotation matrix in the
        $2$-$3$-plane, given by
        \begin{align*}
            C_t = \begin{pmatrix}
                1&0&0&0\\0&\cos(t) & -\sin(t)&0 \\ 0 &\sin(t) &
                \cos(t)&0\\0&0&0&1
            \end{pmatrix}.
        \end{align*}
        We define the universal homotopy $H\colon
        \ures\times [0,\pi/2] \to \ures$, which at time $t$ is
        \begin{align*}
            H_t(X) &= \rho^* C_t^* (X\oplus \flip(X)) C_t \rho
            \\&=\begin{pmatrix}
                X_{++} & \sin(t) X_{-+} & \cos(t) X_{-+} & 0 \\ 
                \sin(t) X_{+-} & X_{--} &0 & \cos(t) X_{+-} \\
                \cos(t) X_{+-} & 0 & X_{--} &-\sin(t) X_{+-} \\
                0 & \cos(t) X_{-+} & -\sin(t) X_{-+} & X_{++}
            \end{pmatrix}.
        \end{align*}
        It is a unitary matrix whose $+-$ and $-+$ components are
        $L^2$-operators, and so it lies in $\ures$. The map
        $H_{\frac{\pi}{2}}$ has values always in the subgroup $\mathrm U_+ \times
        \mathrm U_- \subset \ures$, while $H_{0}$ is just $X\mapsto X\boxplus
        \flip(X)$. This homotopy furthermore induces a well-defined homotopy on
        the quotient
        $\grres$: If we consider another representative $XV = X \begin{pmatrix}
            v_+ & 0 \\ 0 & v_-
        \end{pmatrix}$ for unitary matrices $v_{\pm}$, it is easy to see that
        the operator 
        $V\oplus \flip(V)$ commutes with $C_t$. Therefore, we have 
        \begin{align*}
            \rho^*C_t^* (X V \oplus \flip(X V)) C_t \rho &= \rho^*C_t^* (X \oplus \flip(X))
            (V \oplus \flip (V))C_t\rho \\&= \rho^* C_t^* (X \oplus \flip(X))\,
            C_t \rho \rho^* (V
            \oplus \flip (V)) \rho
        \end{align*}
        and since $\rho^* (V \oplus \flip (V)) \rho \in \mathrm U_+ \times
        \mathrm U_-$, the homotopy is well-defined as a map $\grres\times
        [0,\pi/ 2]\to \grres$. As it goes from $X\boxplus \flip (X)$ to the
        constant map to the basepoint, we are reduced to showing that
        $H$ is a CS-homotopy. The Chern--Simons form is certainly closed, since
        \begin{align*}
            \dd \int_I H_t^* \ch = H_1^*\ch - H_0^*\ch = (X\boxplus \flip X)^*
            \ch - \mathrm{const}_{\h_+}^* \ch = 0
        \end{align*}
        by Prop. \ref{prop:CSforms}. Since $H^{\mathrm{odd}}(\grres)
        = 0$, it follows that it must be exact, and we are done.

        For the odd case, we use the homotopy from Lem. \cite[Lemma
        3.7]{tradler_elementary_2012} in the universal case:
        \begin{align*}
            H_t\colon \mathrm U^1\times I &\to \mathrm U^1\\
            (A,t) &\mapsto (A \boxplus 1) C_t (1\boxplus A^*) C_t^*,
        \end{align*}
        where $C_t = \begin{pmatrix}
            \cos(t) &\sin(t) \\ -\sin(t)&\cos(t)
        \end{pmatrix}$ is again a rotation matrix. This is a homotopy from
        $H_0(A)=A\boxplus A^*$ to $H_{\frac{\pi}{2}}(A)=1$. The calculation
        \begin{align*}
            \dd \int_I H_t^* \ch = H_1^*\ch - H_0^*\ch =
            \mathrm{const}_{1}^* \ch - (A\boxplus A^*)^*
            \ch  = 0
        \end{align*}
        shows that it has a closed CS-form and the vanishing of
        $H^{\mathrm{even}, >0}(\mathrm U^1)=0$ shows that it has to be exact in
        positive degree. In the degree $0$ case, we again make an explicit
        calculation:
        \begin{align*}
            \CS_0(H_t) = \int_I H_t^* \ch_1 =\left(\frac{i}{2\pi}\right)\int_I
            \iota_{\partial_t} ( H_t^* (\tr (\omega_{\mathrm{MC}})))
        \end{align*}
        The integrand resolves to 
        \begin{align*}
            \iota_{\partial_t} &( H_t^* (\tr (\omega_{\mathrm{MC}}))) \\&= \tr
            (C_t (1 \boxplus A)^* C_t^* (A \boxplus 1)^*)((A \boxplus 1)\dot C_t (1
            \boxplus A) X_t^* - (A \boxplus 1)C_t (1 \boxplus A) C_t^* \dot C_t
            C_t^*) \\&= \tr\,C_t^* \dot C_t - \tr \,\dot C_t C_t^* = 0,
        \end{align*}
        thus proving the exactness of the Chern--Simons form.
    \end{proof}
\end{lem}

\section{Natural transformations and exact sequences}
\label{sec:natural}

It remains to define the curvature map $R$, the integration map $I$ and the
action map $a$ from the definition of a differential extension and see that they
have the required properties. Two of those are easy: The underlying class
$I([f])$ of a $\CS$-class $[f]\in \hat K^*(M)$ is given by dropping to the same
equivalence class under the weaker equivalence relation of
homotopy. Note that the stability equivalence relation is also compatible with
homotopy classes by Prop. \ref{prop:spectrum}. The
curvature map $R$ is just the map $\Ch$ from Def. \ref{dfn:chernmaps} which
pulls back the universal Chern form, i.e. $R([f]) = f^*\ch$, where we interpret
$\ch$ as a differential form with values in $V=\R[u,u^{-1}]$, i.e.
\begin{align*}
    \ch_{\mathrm{even}} = \sum_{k\geq 0} \ch_{2k}u^{-k} \in \Omega^0(\grres;V),
    \qquad \ch_{\mathrm{odd}} = \sum_{k\geq 0} \ch_{2k+1} u^{-k}\in \Omega^1(\mathrm
    U^1;V).
\end{align*}
This is well
defined on the $\CS$-equivalence class, since Stokes' theorem implies that
\begin{align*}
    f_1^*\ch - f_0^*\ch = \dd \int_I f_t^*\ch = \dd
    \CS(f_t) = 0,
\end{align*}
if $f_t$ is a $\CS$-homotopy. Stabilization is also fine, since $\Ch(f\boxplus
\mathrm{const}_*) = \Ch(f)$. It is easy to see that 
these maps are homomorphisms. For $R$, this follows from addivity of Chern
forms (Prop. \ref{prop:CSforms}), while for $I$, we use that addition in
$K$-theory can also be implemented
by the blocksum (Prop. \ref{prop:spectrum}). It remains to define the action map 
\begin{align*}
    a\colon \Omega^{*-1}(M;V) / \mathrm{im}(d) \to \hat K^{*}(M).
\end{align*}
Our construction is adapted from \cite[Def. 3.26]{tradler_differential_2015}\cite[Prop.
5.3]{tradler_elementary_2012}. The map has to fit into the exact sequence
\[\begin{tikzcd}
        K^{*-1}(M) \arrow[r, "\mathrm{Ch}"] & \Omega^{*-1}(M;V) / \mathrm{im}(d)
        \arrow[r, "a"] & \hat K^*(M) \arrow[r, "I"]& K^*(M) \arrow[r] &0.
\end{tikzcd}\]
Finding such a map is algebraically equivalent to constructing an isomorphism
\begin{align*}
    \widehat{\CS}: \mathrm{ker}(I) \to \Omega^{*-1}(M;V) / \mathrm{im}(d)
    /\mathrm{im} (\mathrm{Ch}).
\end{align*}
We define $\widehat{\CS}$ as follows: For a map $f_1$ to lie in the kernel of
$I$, it means that there is a homotopy $f_t$ such that
$f_0=\mathrm{const}_{*}$. We define $\widehat \CS(f_1) = \CS(f_t)$ for
a choice of such a nullhomotopy. 
\begin{lem}
    \label{lem:cshat}
    The map $\widehat \CS$ is a well-defined isomorphism of groups, and the
    resulting map
    \[\begin{tikzcd}
        a\colon \Omega^{*-1}(M;V) / \mathrm{im}(d) \ar[r, twoheadrightarrow]&
        \Omega^{*-1}(M;V) / \mathrm{im}(d)/\mathrm{im}(\Ch) \ar[r,
        "\widehat{\CS}^{-1}"] & \mathrm{ker}(I) \ar[r, hookrightarrow] & \hat
        K^{*}(M)
    \end{tikzcd}\]
    meets the axiomatic requirements for the action map in differential $K$-theory.
    \begin{proof}
        Let $g_t$ be another nullhomotopy of $f_1$. We need to show that the
        resulting
        $\CS$-forms only differs by a Chern form up to exact forms. Since both
        homotopies start and end at the same point, we can construct a 
        loop $f_t * g_{1-t}$ at the base point and calculate
        \begin{align*}
            \CS(f_t) - \CS(g_t) = \CS(f_t) + \CS(g_{1-t}) = \CS(f_t *
            {g_{1-t}}).
        \end{align*}
        Recall that we have the homotopy equivalences given by explicit
        periodicity maps 
        \begin{align*}
            \Omega \grres\to \mathrm U^1 \qquad \mathrm{and} \qquad\Omega
            \mathrm U^1\to\grres,
        \end{align*}
        defined in Sec. \ref{sec:smoothper}. We will denote them both
        by the letter $h$. It was shown in Prop.
        \ref{prop:compatibilityeven} and \ref{prop:compatibilityodd} that we
        have 
        \begin{align*}
            \CS(f_t * {g_{1-t}}) = \mathrm{Ch}(h\circ \widehat{(f_t *
                {g_{1-t}})}) + \mathrm{exact}.
        \end{align*}
        Therefore, the map $\widehat \CS$ is well-defined. Furthermore, if we
        have $f_0\boxplus g_0$ with nullhomotopies $f_t$ and $g_t$, then
        \begin{align*}
            \CS(f_t\boxplus g_t) = \CS(f_t) + \CS(g_t)
        \end{align*}
        by Prop. \ref{prop:CSforms} and so $\widehat \CS$ is a homomorphism. 
        
        For injectivity,
        suppose that $\widehat \CS (f_0) = 0$, so $\CS(f_t) = \mathrm{Ch}(g) +
        \mathrm{exact}$ for some $g\colon M\to \mathcal{U}$ for $\mathcal{U}\in
        \left\{ \grres, \mathrm{U}^1 \right\}$. Using any homotopy
        inverse of the periodicity map $h$, we can construct a based loop 
        $\hat{H} := h^{-1}\circ
        g\colon M\to \Omega\mathcal{U}$. Since the domain is just a finite-dimensional
        smooth manifold, we can up to homotopy assume that this map is smooth.
        As an immediate corollary of Prop. \ref{prop:compatibilityeven} and
        \ref{prop:compatibilityodd}, we have that (modulo exact forms)
        \begin{align*}
            \CS(H_t) = \Ch(h\circ \hat H) = \Ch(h \circ h^{-1}\circ g) = \Ch(g).
        \end{align*}
        The composition $H_{t-1} * f_t$ is now a nullhomotopy of $f_0$ with exact
        Chern--Simons form, since
        \begin{align*}
            \CS(H_{1-t} * f_t) = -\CS(H_t) + \CS(f_t) = -\Ch(g) + \Ch(g) +
            \mathrm{exact} = \mathrm{exact}.
        \end{align*}
        We still have to prove surjectivity. The key result here is a
        surjectivity statement for the Chern character on the level of
        differential forms. It is proved in
        \cite[Prop. 2.1 and Rem. 2.3]{Pingali2014} that any exact even form is the Chern form
        of a trivial hermitean bundle with compatible connection. On the other hand, in
        \cite[Cor. 2.7]{tradler_elementary_2012} it is shown that
        that every exact odd form is the Chern form of a nullhomotopic map $M\to\mathrm
        U$\footnote{The proof given there applies without change to the case of
        a noncompact manifold $M$.}. Since our odd Chern character is compatible
        with the \cite{tradler_differential_2015}-Chern character (cf. Sec. \ref{sec:comparison}) and the
        even Chern character is compatible with classifying maps of connections, we
        can conclude that the maps
        \begin{align*}
            \Ch\colon \mathrm{Map}^0(M,\grres) &\to \Omega^0(M;V) \\\Ch\colon
            \mathrm{Map}^0(M,\mathrm U^1) &\to \Omega^1(M;V)
        \end{align*}
        given by pullback of universal forms hit all exact forms. The
        decoration here means that we only consider
        maps which are smoothly homotopic
        to the constant map to the basepoint. The surjectivity of $\CS$ can be
        deduced from this in the following way. Let $\omega$ be an even or odd
        form on $M$. Then, we can construct a form $\widetilde \omega\in
        \Omega(M\times I)$ such that under the inclusions at the endpoints, one
        has $i_0^*\widetilde \omega = 0$ and $i_1^* \widetilde \omega = \omega$.
        The previous result allows us to write $\dd\widetilde \omega =
        \Ch(g_t)$ for some map $g_t$. Via Stokes' theorem, this yields
        \begin{align*}
            \CS(g_t) = \int_I \Ch(g_t) = \int_I \dd \widetilde \omega = \omega +
            \mathrm{exact},
        \end{align*}
        and we are done.
                
        The only thing left to show is that we have
        $R\circ a =d$, but this follows from the calculation
        \begin{align*}
            \dd \int_I f_t^*\ch = f_1^*\ch - f_0^*\ch = f_1^*\ch -
            \mathrm{const}_{\id}^*\ch = f_1^*\ch = R(a(\int_I f_t^*\ch)).
        \end{align*}
    \end{proof}
\end{lem}
With the definitions of $R,a$ and $I$ as above and the abelian group structure
given by blocksum, we have proved our main theorem.
\begin{thm}
    \label{thm:mainthm}
    On the category of possibly non-compact smooth manifolds, the abelian group
    valued functors 
    \begin{align*}
        \hat K^0(M) = \mathrm{Map}(M, \grres) / \CS\mathrm{-homotopy} +
        \mathrm{Stabilization}\\
        \hat K^{1}(M) = \mathrm{Map}(M, \mathrm U^1) / \CS\mathrm{-homotopy} +
        \mathrm{Stabilization}
    \end{align*}
    define differential $K$-theory.
\end{thm}
\begin{rem}
    \label{rem:functoriality}
    Differential $K$-theory is functorial
    for $\CS$-equivalence classes of maps: if two maps $f_0,f_1\colon M\to N$ are
    homotopic and the homotopy $f_t$ satisfies the additional condition that for
    any $g\colon Y\to \grres$ or $g\colon Y\to \mathrm U^1$, the Chern--Simons
    form $\CS(f_t\circ g)$ is exact, then $f_0$ and $f_1$ induce the same map on
    $\hat K^0$ resp. $\hat K^{1}$. This feature of a descent to a quotient category of
    smooth manifolds is a general property of differential cohomology theories
    and is discussed in \cite[Cor. 2.5]{tradler_differential_2015}.
\end{rem}

\section{Comparison to the Tradler-Wilson-Zeinalian model}
\label{sec:comparison}

Though all versions of even differential $K$-theory that meet the
Bunke-Schick axioms are uniquely isomorphic, there are infinitely many
inequivalent versions of odd differential $K$-theory. There is however, only 
one unique isomorphism type that is compatible with the additional structure of an
$S^1$-integration map. One such model for compact manifolds is the one proposed in
\cite[Thm. 4.25]{tradler_differential_2015}. We will establish an explicit isomorphism of both the
even and odd part of the restriction of our model to compact manifolds to the
\cite{tradler_differential_2015}-model in this section.

The \cite{tradler_differential_2015}-model is also based on smooth classifying spaces. For the odd part,
they use the stable unitary group $\mathrm U$ and define 
\begin{align*}
    \hat K_{T}^{1}(M) = \mathrm{Map} (M, \mathrm U) / \CS\mathrm{-equivalence}.
\end{align*}
Since $\mathrm U$ does not admit a
Banach manifold structure, the authors work with universal cocycles given
by the finite-dimensional differential forms (\ref{eq:MCform}) on the
filtration defined by the inclusions of $\mathrm U(n)$ for $n\in \N$. It is
immediately clear that the forms $\ch_{\mathrm{odd}}\in
\Omega^{\mathrm{odd}}(\mathrm U^1)$ that we use pull
back to give the same forms under the natural inclusions
\begin{align*}
    \mathrm U(n)\hookrightarrow \mathrm U \stackrel{i}{\hookrightarrow} \mathrm U^1.
\end{align*}
The second map also preserves the blocksum on the nose. Since $\CS$-homotopies go to
$\CS$-homotopies, it induces a well-defined homomorphism $i_*\colon \hat
K^{1}_{\mathrm{T}}(M)\to  \hat K^{1}(M)$.
\begin{prop}
    \label{prop:comparisoneven}
    The homomorphism $i_*\colon \hat K^1_{T} \to \hat K^1$ preserves all
    the structure of a differential extension, i.e. $I \circ i_* = I_T$, $i_*
    \circ a_T = a$ and $R \circ i_* = R_T$. Furthermore, $i_*$ is an isomorphism.
    \begin{proof}
        The compatibilities are easy to check and follow from $i$ being a
        homotopy equivalence and pulling back $\ch$ to $\ch_T$. The isomorphism
        property of $i_*$ follows from an application of the five lemma to the diagram
        \[\begin{tikzcd}
                K^{0}(M) \arrow[r, "\mathrm{Ch}"]\arrow[d, equal] &
                \Omega^{\mathrm{even}}(M; \R) / \mathrm{im}(d)
                \arrow[r, "a"]\arrow[d, equal] & \hat K_{T}^1(M) \arrow[r,
                "I"]\arrow[d, "i_*"]&
                K^1(M)\arrow[d, equal] \arrow[r] &0\arrow[d, equal] \\
                K^{0}(M) \arrow[r, "\mathrm{Ch}"] & \Omega^{\mathrm{even}}(M; \R) / \mathrm{im}(d)
                \arrow[r, "a"] & \hat K^1(M) \arrow[r, "I"]& K^1(M) \arrow[r] &0.
        \end{tikzcd}\]
    \end{proof}
\end{prop}
The even part of the \cite{tradler_differential_2015}-theory is given by maps into the space of finite rank
projections on $\C^{\infty}_{-\infty} = \bigoplus_\Z \C\subset \h$, defined as
\begin{align*}
    \mathrm{Proj} &= \left\{ \pi\in \mathrm{End}(\C^{\infty}_{-\infty}) \mid \pi^*
    = \pi, \mathrm{Spec}(\pi)\subset \left\{ 0,1 \right\}, \mathrm{rank}(\pi -
    \pi_{\C^{0}_{-\infty}}) < \infty)\right\} 
    \\&\cong  \left\{ V\subset \C^{\infty}_{-\infty}\mid \C^{p}_{-\infty}\subset V
    \subset \C^{q}_{-\infty} \text{ for some } p,q\in\Z\right\}.
\end{align*}
Their basepoint is the space $\C^{0}_{-\infty}$. Apart from a change of basis,
we can identify $\mathrm{Proj}$ with the colimit of the finite-dimensional
Grassmannians, which we denoted by $\mathrm{Gr}_{\mathrm{res},
\mathrm{\infty}}$ in Section \ref{sec:restricted}, as follows: Denote by $A\colon \h\to \h$
the change of basis which maps $e_i$ to $e_{-i}$ for all $i$. Then, we have a natural
map
\begin{align*}
    i\colon \mathrm{Proj} &\to\mathrm{Gr}_{\mathrm{res},
    \mathrm{\infty}}\stackrel{\sim}{\hookrightarrow}
    \grres \\ \pi &\mapsto
    A \pi A
\end{align*}
This is well-defined, since $A\pi A - \pi_+ = A (\pi -
\pi_{\C^{0}_{-\infty}}) A$ has image contained in some $\h_N\subset
\mathrm{im}(A\pi A)\subset \h_{-N}$. We check that it is a homomorphism for the
blocksum. We have that $\pi_1\boxplus \pi_2 = \rho^{*}
\pi_1\oplus\pi_2\rho$ gets mapped to $A \rho^{*} (\pi_1\oplus\pi_2)\rho A $. On
the other hand, the blocksum of the images is $\rho^{*} (A\oplus A) (\pi_1
\oplus \pi_2) (A\oplus A) \rho$. Comparing $(A\oplus A)\rho$ and $\rho A$ as
operators from $\h \to \h\oplus \h$ (see Definition \ref{dfn:blocksum}), we see
that they both map basis vectors $e_{2i}$ to $(e_{-i},0)$. On odd basis vectors,
we have
\begin{align*}
    (A\oplus A)\rho (e_{2i+1}) = (0,e_i), \qquad \rho A (e_{2i+1}) =
    \rho(e_{-2(i+1)+1}) = (0,e_{i+1}).
\end{align*}
Therefore, if we have $f,g\colon M\to
\mathrm{Proj}$, then $i\circ (f\boxplus g)$ and $(i\circ f)\boxplus (i\circ g)$
differ only by conjugation with a fixed unitary matrix $B \in \mathrm U_+\times
\mathrm U_-$ which shifts odd basis vectors by one. By Lem.
\ref{lem:unitaryconjugationisfine}, these are therefore $\CS$-equivalent. We conclude
that $i$ induces a homomorphism of differential $K$-theory groups.
\begin{prop}
    \label{prop:comparisonodd}
    The homomorphism $i_*\colon \hat K^0_{T} \to \hat K^0$ preserves all
    the structure of a differential extension, i.e. $I \circ i_* = I_T$, $i_*
    \circ a_T = a$ and $R \circ i_* = R_T$. Furthermore, $i_*$ is an isomorphism.
    \begin{proof}
        As in the even case, we have that $i$ is a homotopy equivalence. We need
        to check that $i^*\ch = \ch_T$. The path
        components of $\mathrm{Proj}$ are given
        by the rank map, where by definition $\mathrm{rank}(V) =
        \dim(V/\C^{-N}_{-\infty}) - N$, if we and
        assume that $\C^{-N}_{-\infty} \subset V \subset \C^{N}_{-\infty}$. This agrees
        with the path component of the image, which is indexed by $\mathrm{virt.dim}
        (A(V))$. In order to check that also the positive degree parts of
        $\ch$ are compatible, we note that for the inclusion
        $\mathrm{Gr}_{k,2N}\hookrightarrow \mathrm{Proj}$, the \cite{tradler_differential_2015}-Chern character is
        calculated in terms of traces of powers of the differential
        forms $\pi\dd \pi \dd \pi$. Pulling back along the composition 
        \begin{align*}
            \mathrm{Gr}_{k,2N} \to \mathrm{Proj}\to\mathrm{Gr}_{\mathrm{res},
            \mathrm{\infty}} 
        \end{align*}
        on the other hand gives the forms $A \pi \dd \pi \dd \pi A$ whose powers of traces
        agree with those. Since any map from a compact manifold factors through
        one of these Grassmannians, we are done.
    \end{proof}
\end{prop}
We can use the isomorphism $i_*$ in order to slightly improve our model for the
case of a compact manifold: The
stabilization condition that $f\boxplus \mathrm{const}_*$ is equivalent to just $f$ is
actually not needed. We can now prove our second theorem.
\setcounter{thm}{1}
\begin{thm}
    \label{thm:mainresult2}
    On the category of compact smooth manifolds, the abelian group valued functors 
    \begin{align*}
        \hat K^0(M) = \mathrm{Map}(M, \grres) / \CS\mathrm{-homotopy}\\
        \hat K^{1}(M) = \mathrm{Map}(M, \mathrm U^1) /\CS\mathrm{-homotopy}
    \end{align*}
    define differential $K$-theory.
    \begin{proof}
        The only thing left to show is that $f\boxplus \mathrm{const}_* \sim_{\CS} f$ for any
        map $f$. Using the unique isomorphism $(i_*)^{-1}\colon \hat K^* \to
        \hat K_T^*$ to the \cite{tradler_differential_2015}-model, we see that 
        \begin{align*}
            [f\boxplus \mathrm{const}_*] = i_*  ((i_*)^{-1} [f]\boxplus_T
            [\mathrm{const}_*]_T) = i_*(i_*)^{-1}[f] = [f].
        \end{align*}
        In the second equality we used that in the \cite{tradler_differential_2015}-model, we are always on a
        finite stage in the filtration of $\mathrm U$ or $\mathrm{Proj}$. In
        that case, the shuffle blocksum operation as we
        defined it is $\CS$-equivalent to the naive blocksum of
        finite-dimensional matrices, for which the constant map to the identity
        is easily seen to be a unit (compare \cite[Lem. 3.9 and Lem
        3.24]{tradler_differential_2015}).
    \end{proof}
\end{thm}

\section{Examples}
\label{sec:examples}

Already the point is an interesting example, since it illustrates the role that
is played by Chern-Simons homotopies. 
\begin{prop}
    \label{prop:point}
    We have isomorphisms $\hat K^0(*) \cong K^0(*) \cong \Z$ and $\hat
    K^1(*) = \R / \Z$, given by the underlying class map $I$ and the
    determinant map.
    \begin{proof}
        Since there are no odd forms on the point, every homotopy is a
        $\CS$-homotopy, and the first part follows. In the odd case, we need
        to check that the homotopies $f_t$ we use are $\CS$-homotopies, i.e. that the
        $\CS$-form $\CS(f_t) = \int_I f_t^*(\ch_1) = 
        \frac{i}{2\pi} \int_I \tr f_t^{-1}\dot f_t$ is
        exact. We have a splitting induced by the (Fredholm)-determinant map
        \begin{align}
            S\mathrm U^1 \rtimes \mathrm U(1)\cong \mathrm U^1.\label{eq:decomp}
        \end{align}
        Under the isomorphism, the semi-direct group structure is given by
        \begin{align*}
            (n_1, h_1)\cdot(n_2, h_2) = (n_1h_1n_2h_1^{-1}, h_1h_2), \quad (n,h)^{-1} =
            (h^{-1}n^{-1}h,h^{-1}),
        \end{align*}
        which yields that for $f_t = (n_t,h_t)$, we have
        \begin{align*}
            f_t^{-1} \dot f_t = h_t^{-1}n_t^{-1}\dot n_th_t + h_t^{-1}\dot h_t,
        \end{align*}
        where the first term is in $\mathfrak{su}(n)$ and the second one in
        $\mathfrak{u}(1)$. Since $\mathfrak{su}(n)$ consists of matrices with
        trace zero, every homotopy
        that leaves the second factor in (\ref{eq:decomp}) alone will be fine.
        Therefore, all the information is in the second factor, and since the isomorphism
        above is induced by the determinant map, we are done.
    \end{proof}
\end{prop}
Next, we will study the circle $S^1$. One can deduce from the exact sequence
\begin{align*}
    0\to \Omega^{*-1}(S^1;\R)/\mathrm{im}(d)/\mathrm{im}(\Ch) \stackrel{a}{\to}
    \hat K^*(S^1)\stackrel{I}{\to} K^*(S^1)\to 0
\end{align*}
coming from the axioms of $\hat K$ that there are exact sequences
\begin{align}
    0\to H^1(S^1;\R)/\Z \to \hat K^0(S^1) \to K^0(S^1)\to 0\label{eq:exactseq} \\
    0\to \mathscr{C}^\infty(S^1)/\Z \to \hat K^1(S^1) \to K^1(S^1)\to
    0.\nonumber
\end{align}
We will give a description of the kernel of $I$ in both cases. The
following Lemma has already been computed in a different geometric model in
\cite[Lem. 5.3]{bunke2009smooth}.
\begin{lem}
    \label{lem:circle}
    Let $(E_{\pm},\nabla_{\pm})$ be a pair of vector bundles with connection
    over $S^1$ with $\mathrm{dim}(E_+)=\mathrm{dim}(E_-)$. Then, the
    corresponding element in $\hat K^0(S^1)$ is 
    \begin{align*}
        [(E_+,\nabla_+)]-[(E_-,\nabla_-)] = a\left(\frac{1}{2\pi i} \mathrm{log} 
        \frac{\det\mathrm{hol}(E_+,\nabla_+)}{\det\mathrm{hol}(E_-,\nabla_-)}
        \dd z\right),
    \end{align*}
    where $\dd z$ is the volume form on $S^1$ and
    $\mathrm{hol}(E_{\pm},\nabla_{\pm})\in
    \mathrm U(n)/\mathrm{conjugation}$ is
    the holonomy of the bundle. 
    \begin{proof}
        Let $f^\pm$ be the classifying maps of the bundles with connections with
        corresponding nullhomotopies $f_t^\pm$ with $f_0^\pm$ some constant map.
        We get a corresponding homotopy of the classifying map of the virtual
        bundle (cf. Rem. \ref{rem:cyclemap})
        \begin{align*}
            f_t=f_t^{+}\boxplus \mathrm{flip}f_t^{-}\colon S^1\to \grres^0.
        \end{align*}
        By the construction of the map $a$ (Lem. \ref{lem:cshat}), the
        differential form corresponding to $f_1=f$ is the Chern-Simons form of a
        nullhomotopy to $\mathrm{const}_{\h_+}$. It is easy to see that $f_0$,
        which is a constant map to some subspace
        in $\grres^0$, can be connected by a $\CS$-homotopy to
        $\mathrm{const}_{\h_+}$. We are therefore reduced to computing the
        $\CS$-form of $f_t$. Notice that $f_t$ has image contained in
        some finite-dimensional $\mathrm{Gr}_{k,2k}\subset \grres^0$, which
        means that we really have a path of actual finite-dimensional bundles
        with connections. We have
        \begin{align*}
            \omega_f = \int_I f_t^* \ch_2 = \int_I \ch_2(f_t^+) - \ch_2(f_t^-) =
            \frac{i}{2\pi}\int_I \Omega_{\det(f_t^+)} - \Omega_{\det(f_t^-)}.
        \end{align*}
        In the last step, we used that $\ch_2$ of a bundle is the same as the first Chern class of its
        determinant line bundle, i.e. the integral over $\frac{i}{2\pi}$ times its
        curvature. In order to see what is the cohomology class of $\omega_f\in
        \Omega^1(S^1)$, we can just integrate over $S^1$. This yields
        \begin{align*}
            \int_{S^1} \omega_f =
            \frac{i}{2\pi}\int_{S^1}\int_I\Omega_{\det(f_t^+)} -
            \Omega_{\det(f_t^-)} &= \frac{1}{2\pi i} \left(-\int_{D^2}
            \Omega_{\det(f_t^+)}+\int_{D^2}
            \Omega_{\det(f_t^-)}\right) \\&=\frac{1}{2\pi i}
            \left(\mathrm{log}\,
            \mathrm{hol}(\det(f^+))- \mathrm{log}\,
            \mathrm{hol}(\det(f^-))\right),
        \end{align*}
        which shows that the function $f$ corresponds to the logarithm of the
        determinant of the holonomy of its induced bundles, as claimed.
    \end{proof}
\end{lem}
\begin{prop}
    \label{prop:circle2}
    A cohomology class on $S^1$ represented by a one-form $\alpha$ gets mapped in
    the exact sequence (\ref{eq:exactseq}) to $a([\alpha]+\Z[\dd z]) = [f_\alpha]$, where
    $f_\alpha$ is the classifying map of the trivial line bundle with local
    connection form $i \alpha\in \Omega^1(S^1;\mathfrak{u}(1))$. In the odd case, a function $\varphi\in
    \mathscr{C}^{\infty}(S^1)$ gives rise to an
    element in $\hat K^1(S^1)$ via the exponential map, i.e. $a(\varphi + \Z) =
    [\exp(\frac{2\pi}{i} \varphi)]$.
    \begin{proof}
        Let $i\alpha\in \Omega^1(S^1;\mathfrak{u}(1))$ be a local connection
        form for the trivial $\mathrm{U}(1)$-bundle $E$. Furthermore, let 
        \begin{align*}
            s\colon I&\to S^1\times \mathrm{U}(1)\\ t &\mapsto (z,H(t))
        \end{align*}
        be a horizontal lift of the fundamental loop on $S^1$, where
        $z=\exp(2\pi i t)$. If we write
        $\alpha_z =
        \alpha(z) \dd z$, it is determined by the equation
        \begin{align*}
            i\alpha(z)\dd z = H(z)^{-1} H'(z) \dd z,
        \end{align*}
        which we can integrate over the interval and exponentiate in order to get
        \begin{align*}
            \exp\left(i \int_{S^1} \alpha(z)\dd z\right) = \exp\left(\int_{I}
            H(z)^{-1} H'(z) \dd z\right) = H(1).
        \end{align*}
        In the last step, we used that for any path $h\colon I\to \mathrm U(1)$
        starting at $h(0)=1$, we have
        \begin{align*}
            h(s) = \exp\left(\int_0^s h_t^{-1} \dot h_t\dd t\right),
        \end{align*}
        which can be seen by noting that $k(s) = h_s^{-1}\exp(\int_0^s h_t^{-1}
        \dot h_t\dd t)$ satisfies $k(0)=1$ and $\dot k(s)= 0$ for all $s$. Note
        that $H(1)$ is precisely the holonomy of the connection $i\alpha$. By
        Lemma \ref{lem:circle}, we have
        \begin{align*}
            [E,i\alpha] &= a\left(\frac{1}{2\pi i}\mathrm{log}\,\mathrm{H(1)}\,\dd
            z\right) \\&= a\left( \frac{1}{2\pi} \left(\int_{S^1} \alpha(z)\dd
            z\right) \dd z\right),
        \end{align*}
        and the form in the brackets is cohomologous to $\alpha$, which proves
        the claim.

        For the odd part, let $f\colon S^1\to
        \mathrm U^1$ be a representative of a class in $\hat K^1(S^1)$ that is
        in the kernel of $I$. Choose a nullhomotopy $f_t$. The corresponding
        $0$-form $\omega_f$ is the Chern-Simons form of $f_t$, and we calculate
        for $z\in S^1$:
        \begin{align*}
            \omega_f(z) = \int_I (f_t^* \ch_1)(z) = \frac{i}{2\pi} \int_I \tr\,
            f_t^*(z) \dot f_t(z) 
        \end{align*}
        By the splitting in (\ref{eq:decomp}), we can assume that $f_t(z)$
        takes values in $\mathrm U(1)$. Then, for fixed $z$, the integral on the right
        is over a path that starts at $1\in S^1$, and
        we use the same argument as in the even case to conclude
        \begin{align}
            \exp\left(\frac{2\pi}{i} \omega_f(z)\right) = \exp\left(\int_I 
            f_t^*(z) \dot f_t(z)\right) = f(z). \label{eq:exp}
        \end{align}
        Therefore, by defining $f\colon S^1\to \mathrm U^1$ to be the left hand
        side of (\ref{eq:exp}), followed by the inclusion $\mathrm
        U(1)\hookrightarrow
        \mathrm U^1$, we have successfully recovered the function $f$ from the
        given $\omega_f$.
    \end{proof}
\end{prop}

\renewcommand*{\bibfont}{\small}

\printbibliography

\Addresses

\end{document}